\documentclass[11pt]{article}
\usepackage{ulem}
\usepackage{graphicx,amsmath,amsfonts,amssymb,color,mathrsfs, amsmath,amsthm,bm,hyperref}
\usepackage{epsfig}
\newcommand{\chuhao}{\fontsize{19pt}{\baselineskip}\selectfont}

\newcommand{\BOX}{\hfill $\Box$}

\numberwithin{equation}{section}

 \newtheorem{theorem}{Theorem}[section]
 \newtheorem{lemma}{Lemma}[section]
 \newtheorem{assumption}{Assumption}[section]
 \newtheorem{proposition}{Proposition}[section]
 
 \newtheorem{remark}{Remark}[section]

 \newtheorem{example}{Example}[section]
 
 \newtheorem{corollary}{Corollary}[section]

\title{\bf\color{black} \chuhao{Augmented truncation approximations to the solution of Poisson's equation for Markov chains}}

\author{Jinpeng Liu\thanks{School of Mathematics and Statistics, New Campus, Central
South University, Changsha, Hunan, 410083, P.R. China, E-mail: liujinpeng@csu.edu.cn.
\newline \indent \hspace*{-2.5mm}
$^{**}$
School of Mathematics and Statistics, New Campus, Central
South University, Changsha, Hunan, 410083, P.R. China, E-mail:
liuyy@csu.edu.cn.
\newline \indent \hspace*{-2.5mm}
$^{***}$
School of Mathematics and Statistics,
Carleton University, 1125 Colonel By Drive, Ottawa, ON Canada K1S
5B6. Email address: zhao@math.carleton.ca
} \ \ Yuanyuan Liu$^{**}$\ \ Yiqiang Q. Zhao$^{***}$}

\date{Dec  14, 2020}

\begin{document}

\maketitle
\begin{abstract}
Poisson's equation has a lot of applications in various areas.
Usually it is hard to derive the explicit expression of the solution of Poisson's equation for a Markov chain on an infinitely many state space. We will present a computational framework for the solution for both discrete-time Markov chains (DTMCs) and continuous-time Markov chains (CTMCs), by developing the technique of augmented truncation approximations. The convergence to the solution is investigated in terms of the assumption about the monotonicity of the first return times, and is further established for
two types of  truncation approximation schemes: the censored chain and the linear augmented truncation.
Moreover, truncation approximations to the variance constant in central limit theorems (CLTs) are also considered.
The results obtained are applied to discrete-time single-birth processes and continuous-time single-death processes.
\end{abstract}

\section{Introduction}

Let $\mathbf{\Phi}=\{\Phi_k: k\in \mathbb{Z}_+\}$ be a time-homogeneous DTMC on the countable state space $E=\mathbb{Z}_+$
with the probability space $(\Omega,\mathcal{F},{\mathbb{P}})$.
Let $P=(p_{ij})_{i,j\in E}$ be the one-step transition matrix of the chain $\bm{\Phi}$.
Suppose that $P$ is irreducible and positive recurrent with the unique invariant probability vector $\bm{\pi}$
and that $\bm{g}$ is a finite column vector (or function) on $E$ such that $\bm{\pi}^{T}|\bm{g}|<\infty$.
In this paper, we consider truncation approximations to the solution of Poisson's equation and the variance constant in CLTs.
For DTMCs, Poisson's equation has the following form:
\begin{equation}\label{poi-fun-0}
  (P-I)\bm{f}=-\overline{\bm{g}},
\end{equation}
where $I$ is the identity matrix,  $\overline{\bm{g}}=\bm{g}-(\bm{\pi}^{T}\bm{g})\bm{e}$,
and $\bm{e}$ is a column vector of ones.
The vector $\bm{g}$ is called the forcing function, and
the vector $\bm{f}$ is called the solution of Poisson's equation $(\ref{poi-fun-0})$.

Poisson's equation has attracted lots of attention due to its importance in practical applications.
It was pointed out by Meyn and Tweedie \cite{mt09} on pages 458--459  that  Poisson's equation
plays a fundamental role for the analysis of  Markov decision processes, perturbation theory, CTLs, etc.
In Bertsekas \cite{b07}, for a Markov decision process, Poisson's equation was known as the  dynamic programming
equation, and the functions $\bm{g}$ and $\bm{f}$ were called the cost function and the value function, respectively.
In \cite{gm96,l15}, Poisson's equation was adopted for the  perturbation analysis for Markov processes.
In \cite{mt09,sm94}, it was shown that the solution of Poisson's equation can be used to express the variance constant,
which is a very important parameter in CLTs.
Recall that a CLT holds if there exists a constant $0\leq \sigma^{2}(\bm{g}) <\infty$ such that for any initial distribution
\[n^{-\frac{1}{2}}\sum_{k=1}^{n}\overline{g}(\Phi_k)\Rightarrow N(0,\sigma^{2}(\bm{g})),\ \ as \ n\rightarrow \infty,\]
where $N(0,\sigma^{2}(\bm{g}))$ denotes the normal random variable with mean $0$ and variance $\sigma^{2}(\bm{g})$,
and ``$\Rightarrow$''  stands for convergence in distribution.
In addition, Poisson's equation can be  also applied in other fields.
In Glynn and Ormoneit \cite{go02}, Hoeffding's inequality was established for uniformly ergodic DTMCs in terms of the solution of Poisson's equation.
Please refer to \cite{CL19} and references therein for recent developments in this filed.
 Liu and Li \cite{ll18} investigated the error bound for augmented truncation approximations of Markov chains
via Poisson's equation in a discrete-time or continuous-time setting.

For a fixed state $j\in E$,  define
\begin{equation}\label{equ-add-functional}
  f_j ({i})=\mathbb{E}_{i}\bigg[\sum_{k=0}^{\tau_{j}-1}\overline{g}(\Phi_k)\bigg],\ \  i\in E,
\end{equation}
where $\mathbb{E}_{i}[\cdot]:=\mathbb{E}[\cdot|\Phi_0=i]$ denotes the conditional expectation with respect to the initial state $i\in E$
and $\tau_{j}:=\inf\{k\geq 1: \Phi_k= j\}$ is the first return time to state $j$.
It is well known that the vector $\boldsymbol{f_j}$ defined by (\ref{equ-add-functional}) is
a solution of Poisson's equation (\ref{poi-fun-0}) (see, e.g. Glynn and Meyn \cite{gm96}).
For the uniqueness of the solution of Poisson's equation, please refer to Makowski and Shwartz
\cite{MS02} for sufficient criteria.
Since it is not easy to calculate $\boldsymbol{f_j}$ directly, we consider its truncation approximations.

Let $_{(n)}P$ be the $(n+1)\times (n+1)$ northwest
corner truncation of $P$ on $_{(n)}E:=\{0,1,2,\ldots,n\}$.
Let $_{(n)}{\widetilde P}$ be a stochastic transition matrix
such that $_{(n)}{\widetilde P}\geq {_{(n)}P}$
and suppose that $_{(n)}{\widetilde P}$ has a unique invariant probability vector $_{(n)}\bm{\pi}$,
whose corresponding Markov chain is denoted by ${_{(n)}}{\boldsymbol \Phi}=\{ {_{(n)}}\Phi_k: k\in \mathbb{Z}_+\}$.
Let ${_{(n)}}\bm{g}$ be the truncation vector consisting of the first $n+1$ rows of $\bm{g}$.
Similarly, define ${_{(n)}}\tau_{j}$ to be the first return time to state $j$ for  $_{(n)}{\widetilde P}$,
and
\begin{eqnarray*}
_{(n)}f_{j}({i})=\mathbb{E}_{i}\bigg[\sum_{k=0}^{{_{(n)}}\tau_{j}-1}{_{(n)}}\overline{g}({_{(n)}}\Phi_k)\bigg],\ \ i\in {_{(n)}E},
\end{eqnarray*}
where $_{(n)}\overline{\bm{g}}={_{(n)}}\bm{g} - ({_{(n)}}\bm{\pi}^{T}{_{(n)}}\bm{g}){_{(n)}\boldsymbol{e}}$.
The vector $_{(n)}\bm{f_{j}}$ defined above is the unique solution of Poisson's equation
\begin{equation}\label{poi-fun-2}
  ({_{(n)}}\widetilde{P}-{_{(n)}}I){_{(n)}}\bm{f}=-{_{(n)}}\overline{\bm{g}}.
\end{equation}
For truncation approximations, given that $\bm{f_{j}}$ exists, one fundamental issue is to establish the convergence of $_{(n)}\bm{f_{j}}$ to $\bm{f_{j}}$.
Note that $\bm{f_{j}}$ is finite, but it can be unbounded.
We focus on the pointwise convergence, that is, ${_{(n)}}f_{j}(i)\rightarrow f_j(i)$ as $n\rightarrow \infty$ for any $i\in E$.

There exist plenty of literature researches on augmented truncation approximations to invariant probability vectors,
see e.g. recent papers Masuyama \cite{m15,m16} and references therein.
However, to our best knowledge, there is no report on augmented truncation approximations to the solution of Poisson's equation.
The rest of this paper is organized into 5 sections:
In Section 2, we present an example to illustrate that an arbitrarily chosen augmented truncation approximation might not converge to the target solution.
We then investigate the convergence of the augmented truncation approximations in Section 3.
The censored chain and the linear augmentation to some columns are shown to be effective truncation approximation schemes.
We further apply those ideas to approximate the variance constant in CLTs.
In Section 4, we modify the argumentations in Section 3 to  extend the results from DTMCs to CTMCs.
Although most of discrete-time results can be established for the continuous-time case,
we need to pay special attentions to the difference between CTMCs and DTMCs.
In Section 5, we apply our results to single-birth processes and single-death processes,
to derive explicit expressions of the solution of Poisson's equation and the variance constant.
Conclusion and discussion are presented in Section 6.  Two useful propositions are given in Appendix.

\section{An illustrative example}
We consider a DTMC with the following stochastic transition matrix:
\begin{equation*}\label{single-death}
 P=\left (
 \begin{array}{ccccccc}
   q_{0}   & p_{0} & 0     & 0   & \cdots\\
   q_{1}   & 0     & p_{1} & 0   & \cdots\\
   q_{2}   & 0     & 0     & p_2 & \cdots\\
   q_{3}   & 0     & 0     & 0   & \cdots\\
   \vdots&\vdots &\vdots& \vdots & \ddots \\
   \end{array}\right ),
\end{equation*}
where $0<p_i<1$ for each $i\in E$.
Clearly $P$ is irreducible.
Let $p_i=\frac{1}{2}$, if $i=0$ or $i$ is odd, and $p_i=1-\frac{1}{3^{\frac{i}{2}}}$, if $n$ is even.
Define $a_0=1$, $a_i=\prod_{k=0}^{i}p_{k}$, $i\geq1$ and set $\prod_{k=1}^{0}b_{k}=1$.
Then, we have, for any $i\geq1$,
\[
a_i=\left(\frac{1}{2}\right)^{[\frac{i}{2}]+1}\prod_{k=1}^{[\frac{i+1}{2}]-1}\left(1-\frac{1}{3^{k}}\right).
\]
From Liu \cite{l10}, we know that the chain is strongly ergodic.
The invariant probability vector $\bm{\pi}$ is given by
$\pi(0)=1/{\sum_{i=0}^{\infty}a_i}$, $\pi(i)=a_i\pi(0)$, $i\geq1.$
For any finite vector $\bm{g}$ satisfying $\bm{\pi}^{T}|\bm{g}|<\infty$,
there exists a unique solution of Poisson's equation (\ref{poi-fun-0})
due to the special structure of the chain.
Let $j=0$ in (\ref{equ-add-functional}). By calculations, we have
\begin{equation}\label{ill-exam-1}
f_{0}(0)=0,\ \ f_{0}(i)=\left(\bm{\pi}^{T}\bm{g}\right)\sum_{m=0}^{i-1}\frac{1}{\prod_{k=m}^{i-1}{p_{k}}}-
\sum_{m=0}^{i-1}\frac{g(m)}{\prod_{k=m}^{i-1}{p_{k}}},\ \  i\geq 1.
\end{equation}

Now, we consider the last-column augmented matrix $_{(n)}\widetilde{P}$,
which is given by
\begin{equation*}
 _{(n)}\widetilde{P}=\left (
 \begin{array}{ccccccc}
   q_{0}   & p_{0} & \cdots     & 0   & 0\\
   q_{1}   & 0     & \cdots & 0   & 0\\
    \vdots & \vdots     & \ddots   & \vdots & \vdots\\
   q_{n-1}   & 0     & 0     & 0   & p_{n-1}\\
   q_{n}&0 &0& 0 & p_{n} \\
   \end{array}\right ).
\end{equation*}

It is easy to obtain
\[_{(n)}\pi(0)=\frac{1}{\sum_{i=0}^{n-1}a_i+\frac{a_n}{q_n}},\ \
{_{(n)}\pi(i)}=a_i{_{(n)}\pi(0)},\ \  1\leq i\leq n-1,\ \
 _{(n)}\pi(n)=\frac{a_n}{q_n}{_{(n)}}\pi(0),\]
and
\[_{(n)}\bm{\pi}^{T}{_{(n)}\bm{g}}=
\frac{1}{\sum_{i=0}^{n-1}a_i+\frac{a_n}{q_n}}\left(\sum_{i=0}^{n-1}a_ig(i)+\frac{a_n}{q_n}g(n)\right).\]
The unique solution of Poisson's equation for $_{(n)}\widetilde{P}$ is given by
\begin{equation}\label{ill-exam-2}
_{(n)}f_{0}(0)=0,\ \
{_{(n)}}f_{0}(i)=\left({_{(n)}}\bm{\pi}^{T}{_{(n)}}\bm{g}\right)\sum_{m=0}^{i-1}\frac{1}{\prod_{k=m}^{i-1}{p_{k}}}-
\sum_{m=0}^{i-1}\frac{g(m)}{\prod_{k=m}^{i-1}{p_{k}}},\ \ 1\leq i\leq n.
\end{equation}
According to (\ref{ill-exam-1})--(\ref{ill-exam-2}),
the convergence of $_{(n)}\bm{f_{0}}$ to $\bm{f_{0}}$ only depends on the convergence of $_{(n)}\bm{\pi}^{T}{_{(n)}\bm{g}}$
to $\bm{\pi}^{T}\bm{g}$.

We consider the following two choices of the vector $\bm{g}$.
First, let $g(i)=i$ for $i\in E$,
then $\bm{\pi}^{T}|\bm{g}|=\pi(0)\sum_{i=1}^{\infty}ia_i<\infty$.
It implies that the solution $\bm{f_{0}}$ exists.
When $n$ is odd,  $\frac{a_n}{q_n}\rightarrow 0$ and
${_{(n)}}\bm{\pi}^{T}{_{(n)}\bm{g}}\rightarrow\bm{\pi}^{T}\bm{g}$ as $n\rightarrow \infty$;
while when $n$ is even, $\frac{a_n}{q_n}\rightarrow\infty$ and
${_{(n)}}\bm{\pi}^{T}{_{(n)}\bm{g}}\rightarrow\infty$ as $n\rightarrow \infty$.
Thus, we obtain that ${_{(n)}}\bm{\pi}^{T}{_{(n)}\bm{g}}\nrightarrow\bm{\pi}^{T}\bm{g}$ and
${_{(n)}\bm{f_{0}}}\nrightarrow \bm{f_{0}}$ as $n\rightarrow \infty$.
Now, we consider the second choice that $g(i)=i$, if $i$ is odd, and $g(i)=c$, if $i$ is even, where
$c=\sum_{i=0}^{\infty}(2i+1)a_{2i+1}/\sum_{i=0}^{\infty}a_{2i+1}$.
It is easy to verify that
${_{(n)}}\bm{\pi}^{T}{_{(n)}\bm{g}}\rightarrow\bm{\pi}^{T}{\bm{g}}$ as $n\rightarrow\infty$,
from which ${_{(n)}\bm{f_{0}}}\rightarrow \bm{f_{0}}$ as $n\rightarrow \infty$.

\section{Discrete-time Markov chains}

\subsection{General augmented truncations}

In this subsection,
we show the usefulness of the truncation approximations to the solution of Poisson's equation
and to the variance constant for the chain $\boldsymbol{\Phi}$ under the following assumption on the first return times.

Let $j$ be any fixed  state in $E$ and
let $\mathbb{P}_{i}(\cdot):=\mathbb{P}(\cdot|\Phi_0=i)$ be the conditional probability with respect to the initial state $i\in E$.
Define the following additive functionals:
\[
  \zeta_{j}(\bm{g})=\sum_{k=0}^{\tau_{j}-1}{g}(\Phi_k), \ \ {_{(n)}\zeta_{j}({_{(n)}}\bm{g})}=\sum_{k=0}^{{_{(n)}}\tau_{j}-1}{_{(n)}}{g}({_{(n)}}\Phi_k).
\]

\begin{assumption}\label{ass-1}
Suppose that for any initial state $i\in E$ and $n\geq \max \{i, j\}$, both of the following conditions hold:
\begin{description}
  \item[(i)] the sequence $\{_{(n)}\tau_{j}\}$  increases and converges to $\tau_{j}$ with probability one (w.p.1), i.e.
\begin{equation*}\label{tau-1}
\mathbb{P}_{i}\left(\omega\in \Omega: {_{(n)}\tau_{j}}(\omega)\uparrow \tau_{j}(\omega), \ \mbox{as} \  n\rightarrow \infty\right)=1;
\end{equation*}
  \item[(ii)]  the sequence $\{_{(n)}\zeta_{j}(|{_{(n)}\bm{g}}|)\}$  increases and converges to $\zeta_{j}(|\bm{g}|)$ w.p.1, i.e.
\begin{equation*}\label{tau-2}
\mathbb{P}_{i}\left(\omega\in \Omega: {_{(n)}\zeta_{j}(|_{(n)}\bm{g}|)(\omega)}\uparrow \zeta_{j}(|\bm{g}|)(\omega), \ \mbox{as} \  n\rightarrow \infty\right)=1.
\end{equation*}
\end{description}
\end{assumption}

\begin{theorem}\label{the-cor-1}
If Assumption \ref{ass-1} holds, then we have, for any $i\in E$,
\begin{equation}\label{trun-approx-poisson-0}
  \lim_{n\rightarrow\infty}{_{(n)}f_j({i})}=f_j({i}),
\end{equation}
where ${_{(n)}f_{j}({j})}=f_{j}({j})=0$.
\end{theorem}
\proof
For a real number $a$, define
\[a^{+}=\max\{a,0\},\ \ a^{-}=\max\{-a,0\}.\]
Obviously,
$
a=a^{+}-a^{-}.
$
Hence,
\[\mathbb{E}_{i}\left[\zeta_{j}(\bm{g}) \right]
=\mathbb{E}_{i}\left[\zeta_{j}(\bm{g}^{+}) \right]
-\mathbb{E}_{i}\left[\zeta_{j}(\bm{g}^{-})\right].\]

Since $\mathbf{\Phi}$ is positive recurrent, this shows that $\mathbb{E}_i[\tau_{j}]<\infty$ for any $i,j\in E$.
Since $\bm{\pi}^{T}|\bm{g}|<\infty$, from Proposition \ref{tau-con-0} (ii), we know that for any $i, j\in E$,
\begin{equation}\label{zeta-finite}
 |\mathbb{E}_{i}[\zeta_{j}(\bm{g})]| \leq \mathbb{E}_{i}[\zeta_{j}(|\bm{g}|)] <\infty.
\end{equation}
It follows from Theorem 10.31 in \cite{mt09} that
\[
|\bm{\pi}^{T}\bm{g}|=
\left|\frac{1}{\mathbb{E}_j[\tau_{j}]}
\mathbb{E}_{j}\left[\zeta_{j}(\bm{g})\right]\right|
<\infty.
\]
Thus, from Assumption \ref{ass-1} and the monotone convergence theorem, we have
\begin{eqnarray}\label{tau-3-0}
\nonumber \lim_{n\rightarrow\infty}\nonumber {_{(n)}}\bm{{\pi}}^{T}{_{(n)}}\bm{g}
 &=& \lim_{n\rightarrow\infty}\frac{1}{\mathbb{E}_j[{_{(n)}}\tau_{j}]}
    \left( \mathbb{E}_{j}\left[{_{(n)}\zeta_{j}({_{(n)}}\bm{g}^{+})}\right]
     -\mathbb{E}_{j}\left[{_{(n)}\zeta_{j}({_{(n)}}\bm{g}^{-})}\right]\right)\\
\nonumber &=& \frac{1}{\mathbb{E}_j[\lim_{n\rightarrow\infty}{_{(n)}}\tau_{j}]}
  \left(\mathbb{E}_{j}\left[\lim_{n\rightarrow\infty}{_{(n)}\zeta_{j}({_{(n)}}\bm{g}^{+})}\right]-
  \mathbb{E}_{j}\left[\lim_{n\rightarrow\infty}{_{(n)}\zeta_{j}({_{(n)}}\bm{g}^{-})}\right]\right)\\
\nonumber&=&\frac{1}{\mathbb{E}_j[\tau_{j}]}\left(\mathbb{E}_{j}\left[{\zeta_{j}(\bm{g}^{+})}\right]-
  \mathbb{E}_{j}\left[{\zeta_{j}(\bm{g}^{-})}\right]\right)\\
&=&\bm{\pi}^{T}\bm{g}.
\end{eqnarray}
Moreover, we obtain
\[
  f_j (j)=\mathbb{E}_{j}\left[\zeta_{j}(\bm{g})\right] - \left(\bm{\pi}^{T}\bm{g}\right) \mathbb{E}_{j}[\tau_{j}] =0.
\]
Similarly, $_{(n)}f_j (j)=0$ for any $n\geq j$. Hence, we show that (\ref{trun-approx-poisson-0}) holds for $i=j$.

Since $\bm{\pi}^{T}|\bm{g}|<\infty$, from (\ref{zeta-finite}), we know that the solution $f_j({i})$ is finite.
Then, from Assumption \ref{ass-1}, the monotone convergence theorem and (\ref{tau-3-0}), we have, for $i\neq j$,
\begin{eqnarray*}\label{eqn-f-1}
   \lim_{n\rightarrow\infty}{_{(n)}f_{j}({i})}
  &=&\lim_{n\rightarrow\infty}\left(\mathbb{E}_{i}\left[{_{(n)}\zeta_{j}({_{(n)}}\bm{g}^{+})}\right]
     -\mathbb{E}_{i}\left[{_{(n)}\zeta_{j}({_{(n)}}\bm{g}^{-})}\right]
     -\left({_{(n)}}\bm{\pi}^{T}{_{(n)}}\bm{g}\right)\mathbb{E}_{i}\left[{{_{(n)}}\tau_{j}}\right]\right)\\
  &=& \mathbb{E}_{i}\left[\lim_{n\rightarrow\infty}{_{(n)}\zeta_{j}({_{(n)}}\bm{g}^{+})}\right]
   -\mathbb{E}_{i}\left[\lim_{n\rightarrow\infty}{_{(n)}\zeta_{j}({_{(n)}}\bm{g}^{-})}\right]
  -\left(\lim_{n\rightarrow\infty}{_{(n)}}\bm{\pi}^{T}{_{(n)}}\bm{g}\right)\mathbb{E}_{i}\left[\lim_{n\rightarrow\infty}{_{(n)}\tau_{j}}\right]\\
   &=&\mathbb{E}_{i}\left[{\zeta_{j}(\bm{g}^{+})}\right]
   -\mathbb{E}_{i}\left[{\zeta_{j}(\bm{g}^{-})}\right]-(\bm{\pi}^{T}\bm{g})\mathbb{E}_i[\tau_{j}]\\
   &=&f_{j}({i}).
\end{eqnarray*}
So the assertion is proved.
\BOX

From  \cite{mt09}, we immediately know that if $\bm{\pi}^{T}|\bm{g}|<\infty$,
then a CLT holds if for some (then for all, see Proposition \ref{tau-con-0} (i)) $\ell\in{E}$ ,
\begin{equation}\label{var-cos-1-0}
 \mathbb{E}_{\ell}\left[\zeta_{\ell}^2(|\bm{\overline{g}}|)\right]<\infty,
\end{equation}
and if a CLT holds, the variance constant is given by
\begin{equation}\label{var-exp-1-0}
  \sigma^{2}(\bm{g})=\frac{1}{\mathbb{E}_{\ell}[\tau_{\ell}]} \mathbb{E}_{\ell}\left[\zeta_{\ell}^2(\bm{\overline{g}})\right].
\end{equation}

\begin{theorem}\label{variance-trunc}
Suppose that Assumption \ref{ass-1} holds.
If $\mathbb{E}_\ell[\tau_\ell^2]<\infty$ and (\ref{var-cos-1-0}) holds for some $\ell$, then we have
\[
  \lim_{n\rightarrow\infty}{_{(n)}\sigma^{2}({_{(n)}\bm{g}})}=\sigma^{2}(\bm{g}),
\]
where ${_{(n)}\sigma^{2}({_{(n)}\bm{g}})}$ is the variance constant of the chain ${_{(n)}}\boldsymbol{\Phi}$.
\end{theorem}
\proof Using (\ref{var-cos-1-0}) and the following triangle inequality
\[\big||x|-|y|\big|\leq\big|x-y\big|,\]
where $x$ and $y$ are  real numbers, we have
\begin{equation}\label{bound-1-0-0}
\mathbb{E}_{j}\left[\zeta_{j}^2(|\bm{g}|-|\bm{\pi}^{T}\bm{g}|)\right]
\leq \mathbb{E}_{j}\left[\zeta_{j}^2(|\bm{\overline{g}}|)\right]<\infty.
\end{equation}
According to (\ref{bound-1-0-0}) and  the inequality $(x+y)^{2}\leq2(x^{2}+y^{2})$, we obtain
\begin{equation}\label{bound-1-0}
\mathbb{E}_{j}\left[\zeta_{j}^2(|\bm{{g}}|)\right]
\leq 2\mathbb{E}_{j}\left[\zeta_{j}^2(|\bm{g}|-|\bm{\pi}^{T}\bm{g}|)\right]+
2|\bm{\pi}^{T}\bm{g}|^{2}\mathbb{E}_{j}\left[\tau_{j}^2\right]<\infty.
\end{equation}
Using (\ref{bound-1-0}) and H$\ddot{o}$lder inequality derive
\begin{equation}\label{bound-2-0}
\mathbb{E}_{j}\left[\zeta_{j}(|\bm{{g}}|)\cdot\tau_{j}\right] \leq\sqrt{
 \mathbb{E}_{j}\left[\zeta_{j}^{2}(|\bm{{g}}|)\right]\mathbb{E}_{j}\left[\tau_{j}^{2}\right]}<\infty.
\end{equation}

From (\ref{var-exp-1-0}), it is easy to show that for any $n\geq j$,
\[\mathbb{E}_{j}\left[{_{(n)}\zeta_{j}^2({_{(n)}}\bm{\overline{g}})}\right]
=\mathbb{E}_{j}\left[{_{(n)}\zeta_{j}^2({_{(n)}}\bm{{g}})}\right]
-2\left({_{(n)}}\bm{\pi}^{T}{_{(n)}}\bm{g}\right)
\mathbb{E}_{j}\left[{_{(n)}\zeta_{j}({_{(n)}}\bm{{g}})}\cdot{_{(n)}}\tau_{j}\right]
+\left({_{(n)}}\bm{\pi}^{T}{_{(n)}}\bm{g}\right)^{2}\mathbb{E}_{j}[{_{(n)}}\tau_{j}^{2}].\]
Then, by (\ref{bound-1-0})--(\ref{bound-2-0}), Assumption \ref{ass-1} and the monotone convergence theorem, we have
\begin{eqnarray}\label{limit-3-0}
 \nonumber \lim_{n\rightarrow \infty} \mathbb{E}_{j}\left[{_{(n)}\zeta_{j}^2({_{(n)}}\bm{{g}})}\right]
&=& \lim_{n\rightarrow \infty} \mathbb{E}_{j}\left[ {_{(n)}\zeta_{j}^2({_{(n)}}\bm{g}^{+})}-2{_{(n)}\zeta_{j}({_{(n)}}\bm{g}^{+})}{_{(n)}\zeta_{j}({_{(n)}}\bm{g}^{-})}
+{_{(n)}\zeta_{j}^2({_{(n)}}\bm{g}^{-})}\right]\\
\nonumber&=&   \mathbb{E}_{j}\left[{_{(n)}\zeta_{j}^2({_{(n)}}\bm{g}^{+})}-2{\zeta_{j}(\bm{g}^{+})}{\zeta_{j}(\bm{g}^{-})}
+{\zeta_{j}^2(\bm{g}^{-})}\right]\\
&=& \mathbb{E}_{j}\left[{\zeta_{j}^2(\bm{{g}})}\right],
\end{eqnarray}
\begin{equation}\label{limit-2-0}
    \lim_{n\rightarrow \infty} \mathbb{E}_{j}\left[{_{(n)}\zeta_{j}({_{(n)}}\bm{{g}})}\cdot{_{(n)}}\tau_{j}\right] = \mathbb{E}_{j}\left[{\zeta_{j}(\bm{{g}})}\cdot\tau_{j}\right],
\end{equation}
and
\begin{equation}\label{limit-1-0}
 \lim_{n\rightarrow \infty} \mathbb{E}_{j}[{_{(n)}}\tau_{j}^{2}] =\mathbb{E}_{j}[\tau_{j}^{2}],
\end{equation}
where the proof of (\ref{limit-2-0}) is similar to that for (\ref{tau-3-0}).
Hence, by (\ref{tau-3-0}) and (\ref{limit-3-0})--(\ref{limit-1-0}),
\[\lim_{n\rightarrow \infty} {_{(n)}\sigma^{2}({_{(n)}\bm{g}})}
=\lim_{n\rightarrow \infty} \frac{1}{\mathbb{E}_{j}[{_{(n)}}\tau_{j}]} \mathbb{E}_{j}\left[{_{(n)}\zeta_{j}^2({_{(n)}}\bm{\overline{g}})}\right]
=  \frac{1}{\mathbb{E}_{j}[\tau_{j}]}
\mathbb{E}_{j}\left[{\zeta_{j}^2(\bm{\overline{g}})}\right]=\sigma^{2}(\bm{g}).\]
We obtain the assertion of this theorem. \BOX

\begin{remark}\label{the-condition-constant-discrete}
In Theorem \ref{variance-trunc}, we need the extra condition that $\mathbb{E}_\ell[\tau_\ell^2]<\infty$,
which is not very restricted.
Jones and Galin\cite{J04}  showed that (\ref{var-cos-1-0}) holds if $\mathbf{\Phi}$ is geometrically ergodic and $\bm{\pi}^{T}|\bm{g}|^{2+\eta}<\infty$ for some $\eta>0$.
Similar results are given by Roberts and  Rosenthal \cite{rr04}, which showed that  (\ref{var-cos-1-0}) holds if $\bm{\Phi}$ is strongly ergodic and $\bm{\pi}^{T}|\bm{g}|^{2}<\infty$.  Either geometric ergodicity or strong ergodicity implies that $\mathbb{E}_\ell[\tau_\ell^2]<\infty$ for any $\ell\in E$.
\end{remark}

In the following subsections, we consider two special augmented truncations, which will be shown to satisfy Assumption \ref{ass-1}
and can be served as feasible schemes for truncation approximations to the solution of Poisson's equation or to the variance constant.

\subsection{The censored Markov chain}
We first introduce the concept of censoring.
Let $\theta_k$ be the $k$th time that $\Phi_k$  successively visits a state in $_{(n)}E$,
i.e. $\theta_0:=\inf\{m\geq0: \Phi_m\in {_{(n)}E}\}$ and $\theta_{k+1}:=\inf\{m\geq\theta_{k}+1: \Phi_m\in {_{(n)}E}\}$.
The censored Markov chain ${_{(n)}}{\boldsymbol \Phi}=\{ {_{(n)}}\Phi_k: k\in \mathbb{Z}_+\}$ on  $_{(n)}E$
is defined by ${_{(n)}}\Phi_k=\Phi_{\theta_{k}}$, $k\in \mathbb{Z}_+$.
Let $P_{_{E_1,E_2}}=(p_{ij})_{i\in E_1,j\in E_2}$, where  $E_1$ and $E_2$ are subsets of $E$.
According to $_{(n)}E$ and its complement $_{(n)}E^{C}$, we partition the transition matrix $P$ as
\begin{equation*}
P=
 \left(\aligned
 &P_{_{(n)}E, {_{(n)}E}}\ &&P_{{_{(n)}E}, {_{(n)}E^{C}}}\\
 &P_{{_{(n)}E^{C}}, {_{(n)}E}}\ && P_{_{(n)}E^{C}, {_{(n)}E^{C}}} \endaligned \right).
\end{equation*}
Note that $P_{_{(n)}E,{_{(n)}E}}={_{(n)}P}$.
The transition matrix of ${_{(n)}}{\boldsymbol \Phi}$ is  given by (see e.g. page 118 of Latouche and Ramaswami \cite{gv99}):
\begin{equation}\label{gen-tru-2-0}
_{(n)}\widetilde{P}={_{(n)}P}+P_{{_{(n)}E},{_{(n)}E^{C}}}(I-P_{_{(n)}E^{C},{_{(n)}E^{C}}})^{-1}P_{{_{(n)}E^{C}},{_{(n)}E}},
\end{equation}
where $(I-P_{_{(n)}E^{C},{_{(n)}E^{C}}})^{-1}P_{{_{(n)}E^{C}},{_{(n)}E}}$ denotes the probability of first hitting ${_{(n)}E}$ from ${_{(n)}E^{C}}$.

\begin{lemma}\label{lem-con-1-0}
Assumption \ref{ass-1} holds for the censored Markov chain $_{(n)}\widetilde{P}$ defined by (\ref{gen-tru-2-0}).
\end{lemma}
\proof
First we show that for any initial state $i$, if the state $j$ can be reached for $\bm{\Phi}$,
then the state $j$ must be first reached for $_{(n)}\bm{\Phi}$ in a shorter time.
To show this, we adopt the definition of local time on page 118 of \cite{gv99}.
We assume that there are two clocks in the system. The first one measures the time of the chain $\bm{\Phi}$
on the whole state space $E$ and is called the global clock.
The second clock, called the local clock, measures the time in the set $_{(n)}E$.
The local clock increases by one unit per unit of global clock during those intervals
when $\bm{\Phi}$ is in $_{(n)}E$ and remains constant when $\bm{\Phi}$  is  in $_{(n)}E^{C}$.
As is shown in  Figure 1, for a sample path $\omega=(\omega_k, k\in \mathbb{Z}_+)$,
the time intervals of $_{(n)}\bm{\Phi}$, $_{(n+1)}\bm{\Phi}$ and $\bm{\Phi}$, which start at $i$ and end at $j$,
are respectively made up of the cross-shaped points, the cross-shaped and the triangle-shaped points, and all the points.
Thus, for each $\omega$, $_{(n)}\tau_{j}(\omega)$ is an increasing sequence.
Since $\bm{\Phi}$ is irreducible and  recurrent, this shows that $\mathbb{P}_i(\omega\in \Omega: \tau_j(\omega)<\infty)=1$ for any pair of states $i$ and $j$.
Hence we have that for any initial state $i$ and for  each $\omega$,
\begin{equation}\label{gen-tru-3-0}
  _{(n)}{\tau}_{j}(\omega)\uparrow\tau_j(\omega)<\infty,\ \ \mbox{as} \ n\rightarrow \infty.
\end{equation}
Thus, we have verified (i) of Assumption \ref{ass-1}.

\begin{figure}[h]\label{1-1-1}
  \centering
  \includegraphics[width=8cm]{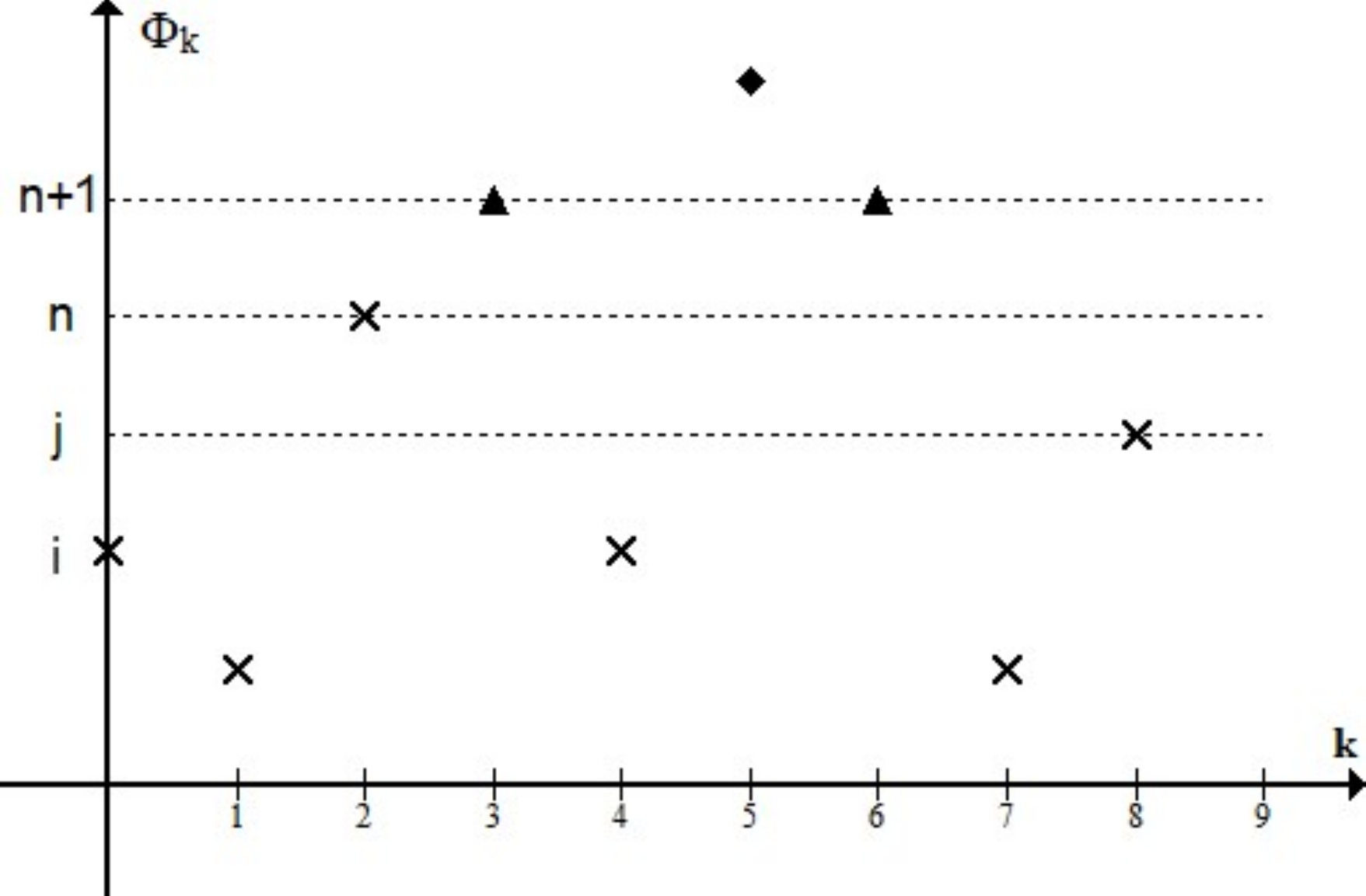}\\
  \caption{Illustrate the first return time to $j$  for $_{(n)}\bm{\Phi}$, $_{(n+1)}\bm{\Phi}$ and $\bm{\Phi}$ through one sample path,
given that the initial state is $i$ and $_{(n)}\widetilde{P}$ is defined by (\ref{gen-tru-2-0}).}
\end{figure}

Now, we verify (ii) of Assumption \ref{ass-1}. Since $\bm{\Phi}$ and $_{(n)}\bm{\Phi}$ are irreducible and recurrent, for any pair of states $i$ and $j$,
\[\mathbb{P}_i(\omega\in \Omega: \tau_{j}(\omega)<\infty)=1\ \ \mbox{and}\ \ \mathbb{P}_i(\omega\in \Omega: {_{(n)}\tau_{j}}(\omega)<\infty)=1.\]
Hence, it immediately follows that
\begin{equation*}\label{pro-1-ii}
  \mathbb{P}_i(\omega\in \Omega: \{\tau_{j}(\omega)<\infty\}\cap\{{_{(n)}\tau_{j}}(\omega)<\infty\})=1.
\end{equation*}
Let $A=\{\omega\in \Omega:\{\tau_{j}(\omega)<\infty\}\cap\{{_{(n)}\tau_{j}}(\omega)<\infty\}\}$,
and note that $\mathbb{P}_{i}(A)=1$. For any $\omega\in A$,
since the vector $\bm{g}$ is finite,
both ${_{(n)}\zeta_{j}(|{_{(n)}\bm{g}}|,\omega)}$ and $\zeta_{j}(|\bm{g}|,\omega)$ exist simultaneously and
they are finite for any $n\geq \max \{i, j\}$.
Similar to the analysis of (\ref{gen-tru-3-0}), we have
\begin{equation*}\label{gen-tru-45}
  {_{(n)}\zeta_{j}(|{_{(n)}\bm{g}}|,\omega)}\leq{_{(n+1)}\zeta_{j}(|{_{(n+1)}\bm{g}}|,\omega)}\leq \zeta_{j}(|\bm{g}|,\omega).
\end{equation*}
Let $M=\max\{|g({\Phi}_{k}(\omega))|: 0\leq k\leq \tau_{j}(\omega)\}$.
From the first assertion, for any $\varepsilon>0$, there exists some $N\in \mathbb{Z}_{+}$, such that for $n>N$,
\[
0<\tau_{j}(\omega)-{_{(n)}\tau_{j}(\omega)}<\varepsilon/M.
\]
Moreover, we obtain
\[
|\zeta_{j}(|\bm{g}|,\omega)-{_{(n)}\zeta_{j}(|{_{(n)}\bm{g}}|,\omega)}|
\leq\sum_{k={_{(n)}}\tau_{j}(\omega)}^{\tau_{j}(\omega)}M
< M \frac{\varepsilon}{M} = \varepsilon.
\]
Then the proof is completed. \BOX

\subsection{Linearly augmented truncation}

We now consider the  linear augmentation, suggested in Seneta \cite{s80},
by augmenting the truncated transition elements to some column.
Specifically, for the fixed state $j$, the transition matrix of the $(j+1)$th column augmented Markov chain is given by
\begin{equation}\label{gen-tru-2-1}
_{(n)}\widetilde{P} ={_{(n)}P}+({_{(n)}I}-{_{(n)}P}){_{(n)}\boldsymbol{e}}{_{(n)}\boldsymbol{e}_{j}}^{T}, \ \ n \geq j,
\end{equation}
where $_{(n)}\boldsymbol{e}_{j}$ is a $(n+1)$-vector with unity in the $(j+1)$th position, zeros elsewhere.
By \cite{s80}, we know that $_{(n)}\widetilde{P}$ has only one closed class in the state space $_{(n)}E$.
It implies that $_{(n)}\widetilde{P}$ has the unique invariant probability vector ${_{(n)}}\bm{{\pi}}$.
\begin{lemma}\label{lem-2}
Assumption \ref{ass-1} holds for the linearly augmented  Markov chain $_{(n)}\widetilde{P}$ defined by (\ref{gen-tru-2-1}).
\end{lemma}
\proof
Similar to the proof of  Lemma \ref{lem-con-1-0}, our arguments are also based on the sample path analysis.
For any sample path $\omega=(\omega_k, k\in \mathbb{Z}_+)$ and initial state $i$,
if the state $j$ can be reached without going through any state in $_{(n)}E^{C}$,
then $_{(n)}\tau_{j}(\omega)$ is equal to $\tau_j(\omega)$.
Otherwise, $_{(n)}\tau_{j}(\omega)$ is small than $\tau_j(\omega)$.
Since the time $_{(n)}\tau_{j}(\omega)$ to reach $j$ in the chain {$_{(n)}\bm{\Phi}$ is calculated by collapsing all time segments
when the chain arrives outside $_{(n)}E$. Using the same argument and replacing $\bm{\Phi}$ by $_{(n+1)}\bm{\Phi}$}, we conclude
that $_{(n)}\tau_{j}(\omega)$ is increasing (See Figure 1 for a depiction).
The rest proof  is similar to that for  Lemma \ref{lem-con-1-0}, which is omitted here.
\BOX

\section{Continuous-time Markov chains}

We now modify the augmentations in Section 3 to adapt the analysis from DTMCs to CTMCs.
Let $\mathbf{\Phi}=\{\Phi_t: t\in \mathbb{R}_+\}$ be an irreducible and time-homogeneous CTMC on a countable
state space $E$ with the $q$-matrix $Q=(q_{ik})_{i,k\in E}$ and the transition function $P_{ik}(t)$.
Assume that $Q$ is totally stable and regular.
We further assume that $\bm{\Phi}$ is positive recurrent with the unique invariant probability vector $\bm{\pi}$.
For CTMCs, Poisson's equation has the following form:
\begin{equation}\label{poi-fun-1}
  Q\bm{f}=-\overline{\bm{g}},
\end{equation}
and we say a CLT holds if there exists a constant $0\leq \sigma^{2}(\bm{g}) <\infty$
such that for any initial distribution
\[t^{-\frac{1}{2}}\int_{0}^{t}\overline{g}(\Phi_s)ds\Rightarrow N(0,\sigma^{2}(\bm{g})),\ \ as \ n\rightarrow \infty.\]

Let $j$ be any fixed state in $E$. We define
\[\xi_{j}(\bm{g})=\int_{0}^{\delta_{j}}g(\Phi_t)dt,\ \
f_j ({i})=\mathbb{E}_{i}\left[\xi_{j}(\bm{\overline{g}})\right],\ \ i\in E,\]
where $\delta_{j}:=\inf\{t\geq J_1: \Phi_t=j\}$ is the first return time to state $j$ and $J_1$ is the first jump time of $\bm{\Phi}$.
According to Asmussen and Bladt \cite{sm94}, we know that the vector $\bm{f_{j}}$ is a solution of Poisson's equation (\ref{poi-fun-1}).

Let $_{(n)}Q$ be the $(n+1)\times (n+1)$ northwest corner truncation of $Q$  on $_{(n)}E$.
Denote by $_{(n)}{\widetilde Q}$ a conservative
$q$-matrix such that $_{(n)}{\widetilde Q}\geq {_{(n)}Q}$ and  $_{(n)}{\widetilde Q}$
has a unique invariant probability vector $_{(n)}\bm{\pi}$,
whose corresponding Markov chain is denoted by ${_{(n)}}\mathbf{\Phi}=\{{_{(n)}}\Phi_t: t\in \mathbb{R}_+\}$.
Let ${_{(n)}J_1}$ and ${_{(n)}\delta_{j}}$ be the first jump time and the first return time to state $j$
of ${_{(n)}}\bm{\Phi}$, respectively. Define
\[{_{(n)}\xi_{j}({_{(n)}}\bm{g})}=\int_{0}^{{_{(n)}\delta_{j}}}{_{(n)}}g(_{(n)}{\Phi}_t)dt,\ \
_{(n)}f_{j}({i})= \mathbb{E}_{i}\left[{_{(n)}\xi_{j}({_{(n)}}\bm{\overline{g}})}\right], \ \ i\in {_{(n)}E},\]
where the vector  $_{(n)}\bm{f_{j}}$ is the unique solution of Poisson's equation:
\begin{equation}\label{poi-fun-3}
  {_{(n)}}\widetilde{Q} {_{(n)}}\bm{f} =- {_{(n)}}\overline{\bm{g}}.
\end{equation}

In the following subsections, we consider the convergence of ${_{(n)}}\bm{f_{j}}$ to $\bm{f_{j}}$.
Due to the difference, see Remark \ref{the-difference}, between DTMCs and CTMCs,
we cannot establish a framework which unifies the treatment for both censored Markov chain and the linearly augmented truncation.
Hence, we consider the two cases separately.

\subsection{The censored Markov chain}

We now introduce the concept of censoring for a CTMC.
To define the censored Markov chain ${_{(n)}}\mathbf{\Phi}$ on ${_{(n)}}E$,
we first define the local clock which increases by one unit per unit global time during those intervals
when $\Phi_t$ is in ${_{(n)}}E$ and remains constant when $\Phi_t$ is not in ${_{(n)}}E$.
Then, we define ${_{(n)}}\Phi_t=\Phi_{t'}$ if $t'$ is the global time corresponding to the local time $t$.
Let $Q_{_{E_1,E_2}}=(q_{ij})_{i\in E_1,j\in E_2}$. Then, we partition the $q$-matrix as
\begin{equation*}
Q=
 \left(\aligned
 &Q_{_{(n)}E,{_{(n)}E}}\ &&Q_{{_{(n)}E},{_{(n)}E^{C}}}\\
 &Q_{{_{(n)}E^{C}},{_{(n)}E}}\ && Q_{_{(n)}E^{C},{_{(n)}E^{C}}} \endaligned \right),
\end{equation*}
where $Q_{_{(n)}E,{_{(n)}E}}={_{(n)}Q}$. According to page 126 of \cite{gv99}, the generator of ${_{(n)}}\bm{\Phi}$ is given by
\begin{equation}\label{gen-tru-2}
_{(n)}\widetilde{Q}={{_{(n)}Q}}+Q_{{_{(n)}E},{_{(n)}E^{C}}}(-Q_{_{(n)}E^{C},{_{(n)}E^{C}}})^{-1}Q_{{_{(n)}E^{C}},{_{(n)}E}}.
\end{equation}
It is known from Theorem 5.5.3 in \cite{gv99} that the censored chain ${_{(n)}}\bm{\Phi}$ is irreducible and positive recurrent with the
unique invariant probability vector $_{(n)}\bm{\pi}$, which is given by
\begin{equation}\label{sencor-pi-1}
  {_{(n)}\pi(i)}=\frac{\pi(i)}{\sum_{k=0}^{n}\pi(k)},\ \ i\in {_{(n)}E}.
\end{equation}

\begin{theorem}\label{the-cor-3-1}
Let $_{(n)}\widetilde{Q}$ be defined by (\ref{gen-tru-2}). Then
\begin{description}
  \item[(i)] for any initial states $i\neq j\in E$ and $n\geq \max \{i, j\}$,
both of the sequences $\{_{(n)}\delta_{j}\}$ and $\{_{(n)}\xi_{j}(|{_{(n)}}\bm{g}|)\}$
increase and converge to $\delta_{j}$ and $\xi_{j}(|\bm{g}|)$ w.p.1, respectively; and
  \item[(ii)] for any states $i\in E$,
\begin{equation}\label{trun-approx-poisson}
  \lim_{n\rightarrow+\infty}{_{(n)}f_j({i})}=f_j({i}).
\end{equation}
\end{description}
\end{theorem}

\proof
The proof of the first assertion is similar to that for Lemma \ref{lem-con-1-0} for the case $i\neq j$,
which is omitted here.

Now, we consider the second assertion. Since $\bm{\pi}^{T}|\bm{g}|<\infty$, it follows from (\ref{sencor-pi-1}) that
\begin{equation}\label{tau-3}
  \lim_{n\rightarrow\infty}{_{(n)}}\bm{\pi}^{T}{_{(n)}}\bm{g}=\lim_{n\rightarrow\infty}\frac{\sum_{i=0}^{n}\pi(i)g(i)}
  {\sum_{k=0}^{n}\pi(k)}=\sum_{i=0}^{\infty}\pi(i)g(i)=\bm{\pi}^{T}\bm{g}.
\end{equation}
From Theorem 1.2 in Section 6.1 of Asmussen \cite{a03}, we have
\[
  f_j (j)=\mathbb{E}_{j}\left[\xi_{j}(\bm{g})\right] - \left(\bm{\pi}^{T}\bm{g}\right)\mathbb{E}_{j}[\delta_{j}]=0.
\]
Similarly, $_{(n)}f_j (j)=0$ for all $n\geq j$. Hence, we showed that (\ref{trun-approx-poisson}) holds for $i=j$.

If  $i\neq j$,  by the first assertion and (\ref{tau-3}),
the proof is analogous to that for Theorem \ref{the-cor-1}.
Thus, this completes the proof.
\BOX

\begin{remark}\label{the-difference}
We cannot expect that (i) of Theorem \ref{the-cor-3-1} holds for the case of $i=j$, which is different from the discrete-time case.
We now explain  why we have to be careful for this situation.
Consider a $q$-matrix $Q$ (see Liu et al. \cite{lzz10})  such that
$q_{00}=-\lambda_{0}$, $q_{0i}=\lambda_{0}p_i$, $q_{0i}=-q_{ii}=\lambda_{i}$ for $i\geq1$, and $q_{ij}=0$ for other $i, j\in E$,
where $\sum_{i\geq1}p_i=1$. Suppose that $Q$ is positive recurrent, i.e. $\sum_{i\geq1}p_i\lambda_{i}^{-1}<\infty$.
According to (\ref{gen-tru-2}), we have, for any $n\in E$,
\[
_{(n)}\widetilde{Q}=\left (
 \begin{array}{ccccc}
  -\lambda_{0}\sum_{i=1}^{n}p_{i} & \lambda_{0}p_{1}&\lambda_{0}p_{2} & \cdots &\lambda_{0}p_{n}  \\
  \lambda_{1} & -\lambda_{1} &0 & \cdots &0 \\
  \lambda_{2} &0 & -\lambda_{2}  & \cdots &0 \\
  \vdots&\vdots & \vdots & \ddots & \vdots\\
  \lambda_{n}&0&0&0&-\lambda_{n}\\
   \end{array}\right ).
 \]
By simple calculations, we have
\[
  \mathbb{E}_{0}[{_{(n)}}\delta_{0}]=\frac{1}{\sum_{i=1}^{n}p_{i}}\bigg(\frac{1}{\lambda_{0}}+\sum_{i=1}^{n}\frac{p_{i}}{\lambda_{i}}\bigg).
\]
Now, let $\lambda_{i}=2^{i}$ for $i\in E$, and $p_i=\frac{1}{2^{i}}$ for $i\geq1$.
Then $\mathbb{E}_{0}[{_{(n)}}\delta_{0}]=\frac{4}{3}(1-\frac{1}{4^{n+1}})/(1-\frac{1}{2^{n}})$,
which  is decreasing in $n$. It implies that the sequences $\{_{(n)}\delta_{0}\}$ cannot increase w.p.1 for the case $i=j=0$.
\end{remark}

From  Glynn and Whitt \cite{Glyn2002}, we immediately know that if $\bm{\pi}^{T}|\bm{g}|<\infty$, then a CLT holds if
 for some (then for all, see Proposition \ref{tau-con-1} (i)) $\ell\in{E}$
\begin{equation}\label{var-cos-1}
 \mathbb{E}_{\ell}\left[\xi_{\ell}^{2}(|\bm{\overline{g}}|)\right]<\infty.
\end{equation}
From \cite{sm94}, if a CLT holds, the variance constant is given by
\begin{equation}\label{var-exp}
\sigma^{2}(\bm{g})=\frac{1}{\mathbb{E}_{\ell}[\delta_{\ell}]}
\mathbb{E}_{\ell}\left[\xi_{\ell}^{2}(\bm{\overline{g}})\right]
=2\sum_{i=0}^{\infty}\pi(i)\overline{g}(i){f_{\ell}(i)}.
\end{equation}

\begin{theorem}\label{variance-trunc-2}
Let $_{(n)}\widetilde{Q}$ be defined by (\ref{gen-tru-2}).
If $\mathbb{E}_\ell[\delta_\ell^2]<\infty$ and (\ref{var-cos-1}) holds for some $\ell\in E$, then we have
\[\lim_{n\rightarrow+\infty}{_{(n)}{\sigma}^{2}({_{(n)}\bm{g})}}=\sigma^{2}(\bm{g}).\]
\end{theorem}
\proof Similar to (\ref{bound-1-0-0}) and (\ref{bound-1-0}), we have, for any  $j\in{E}$,
\begin{equation}\label{var-cos-2}
 \mathbb{E}_{j}\left[\xi_{j}^{2}(|\bm{{g}}|-|\bm{\pi}^{T}\bm{g}|)\right]<\infty\ \ \mbox{and}\ \
 \mathbb{E}_{j}\left[\xi_{j}^{2}(|\bm{{g}}|)\right]<\infty.
\end{equation}
Form (\ref{var-cos-2}) and Lemma 3.1 (i) in \cite{sm94}, we have
\begin{equation}\label{var-cos-3}
  \sum_{i=0}^{\infty}\pi(i)(|g(i)|-|\bm{\pi}^{T}\bm{g}|)
\mathbb{E}_{i}\left[\xi_{j}(|\bm{{g}}|-|\bm{\pi}^{T}\bm{g}|)\right]<\infty,
\end{equation}
and
\begin{equation}\label{the-bound-of-var-1}
  \sum_{i=0}^{\infty}\pi(i)|g(i)|\mathbb{E}_{i}\left[\xi_{j}(|\bm{{g}}|)\right]<\infty.
\end{equation}
Since $\mathbb{E}_j[\delta_j^2]<\infty$,  from Theorem 2.1 in \cite{lzz10}, we know that
\begin{equation}\label{the-bound-of-var-2}
  \sum_{i=0}^{\infty}\pi(i)\mathbb{E}_{i}\left[\delta_{j}\right]<\infty.
\end{equation}
Expanding the left hand side of (\ref{var-cos-3}) and using (\ref{the-bound-of-var-1})--(\ref{the-bound-of-var-2}), we thus obtain
\begin{equation}\label{the-bound-of-var-3}
  \sum_{i=0}^{\infty}\pi(i)|g(i)|\mathbb{E}_{i}\left[\delta_{j}\right]<\infty \ \ \mbox{and}\ \
  \sum_{i=0}^{\infty}\pi(i)\mathbb{E}_{i}\left[\xi_{j}(|\bm{{g}}|)\right]<\infty.
\end{equation}

By simple calculations, we have, for any $n\geq j$,
\begin{eqnarray*}
  _{(n)}{\sigma}^{2}({_{(n)}\bm{g}}) &=& 2\sum_{i=0}^{n}{_{(n)}}{\pi}(i)g(i)\mathbb{E}_{i}\left[{_{(n)}\xi_{j}({_{(n)}}\bm{g})}\right]
  -2\left({_{(n)}}\bm{\pi}^{T}{_{(n)}}\bm{g}\right)\sum_{i=0}^{n}{_{(n)}}{\pi}(i)g(i)\mathbb{E}_{i}\left[{_{(n)}}{\delta}_{j}\right]\\
  &&-2\left({_{(n)}}\bm{\pi}^{T}{_{(n)}}\bm{g}\right)
  \sum_{i=0}^{n}{_{(n)}}{\pi}(i)\mathbb{E}_{i}\left[{_{(n)}\xi_{j}({_{(n)}}\bm{g})}\right]
  +2\left({_{(n)}}\bm{\pi}^{T}{_{(n)}}\bm{g}\right)^{2}\sum_{i=0}^{n}{_{(n)}}{\pi}(i)\mathbb{E}_{i}\left[{_{(n)}}{\delta}_{j}\right].
\end{eqnarray*}
By (\ref{sencor-pi-1}), (\ref{the-bound-of-var-1}), Theorem \ref{the-cor-3-1} and the monotone convergence theorem, we obtain
\begin{eqnarray}\label{the-var-of-censor-1}
\nonumber  &&\lim_{n\rightarrow\infty}\sum_{i=0}^{n}{_{(n)}}{\pi}(i)g(i)
\mathbb{E}_{i}\left[{_{(n)}\xi_{j}({_{(n)}}\bm{g})}\right]\\
\nonumber &=& \lim_{n\rightarrow\infty}\frac{1}{\sum_{k=0}^{n}\pi(k)}
\sum_{i=0}^{\infty}\pi(i)g(i)\left(\mathbb{E}_{i}\left[{_{(n)}\xi_{j}({_{(n)}}\bm{g}^{+})}\right]
-\mathbb{E}_{i}\left[{_{(n)}\xi_{j}({_{(n)}}\bm{g}^{-})}\right]\right)\mathbb{I}_{\{i\leq n\}}\\
\nonumber&=&\sum_{i=0}^{\infty}\pi(i)g(i)\left(\mathbb{E}_{i}\left[{\xi_{j}(\bm{g}^{+})}\right]
-\mathbb{E}_{i}\left[{\xi_{j}(\bm{g}^{+})}\right]\right)\\
&=&\sum_{i=0}^{\infty}\pi(i)g(i)\mathbb{E}_{i}\left[{\xi_{j}(\bm{g})}\right].
\end{eqnarray}
Similarly, it follows from (\ref{the-bound-of-var-2})--(\ref{the-bound-of-var-3}) that
\begin{equation}\label{the-var-of-censor-2}
  \lim_{n\rightarrow\infty}\left({_{(n)}}\bm{\pi}^{T}{_{(n)}}\bm{g}\right)
  \sum_{i=0}^{n}{_{(n)}}{\pi}(i)g(i)\mathbb{E}_{i}\left[{_{(n)}}{\delta}_{j}\right]
  =\left(\bm{\pi}^{T}\bm{g}\right)\sum_{i=0}^{\infty}\pi(i)g(i)\mathbb{E}_{i}\left[\delta_{j}\right],
\end{equation}
\begin{equation}\label{the-var-of-censor-3}
  \lim_{n\rightarrow\infty}\left({_{(n)}}\bm{\pi}^{T}{_{(n)}}\bm{g}\right)\sum_{i=0}^{n}{_{(n)}}{\pi}(i)
  \mathbb{E}_{i}\left[{_{(n)}\xi_{j}({_{(n)}}\bm{g})}\right]
  =\left(\bm{\pi}^{T}\bm{g}\right)\sum_{i=0}^{\infty}\pi(i)\mathbb{E}_{i}\left[{\xi_{j}(\bm{g})}\right]
\end{equation}
and
\begin{equation}\label{the-var-of-censor-4}
  \lim_{n\rightarrow\infty}\left({_{(n)}}\bm{\pi}^{T}{_{(n)}}\bm{g}\right)^{2}
  \sum_{i=0}^{n}{_{(n)}}{\pi}(i)\mathbb{E}_{i}\left[{_{(n)}}{\delta}_{j}\right]
  =\left(\bm{\pi}^{T}\bm{g}\right)^{2}\sum_{i=0}^{\infty}\pi(i)\mathbb{E}_{i}\left[\delta_{j}\right].
\end{equation}
By the above limits (\ref{the-var-of-censor-1})--(\ref{the-var-of-censor-4}), we obtain the assertion.
\BOX

\begin{remark}\label{the-condition-constant}
Theorem 3.1 in Liu and Zhang \cite{lz15} showed that (\ref{var-cos-1}) holds if
 $\bm{\Phi}$ is exponentially ergodic and $\bm{\pi}|\bm{g}|^{2+\eta}<\infty$ for some $\eta>0$ or
 $\bm{\Phi}$ is strongly ergodic and $\bm{\pi}|\bm{g}|^{2}<\infty$. Hence, we have similar remarks to that in Remark \ref{the-condition-constant-discrete}.
\end{remark}

\subsection{Linearly augmented truncation}

In this subsection, we consider the linear augmentation for a CTMC.
Similar to DTMCs, for the fixed state $j$,
the $(j+1)$th column augmentation $q$-matrix $_{(n)}\widetilde{Q}$ is given by
\begin{equation}\label{aug-tru-2}
_{(n)}\widetilde{Q} ={_{(n)}Q}+(-{_{(n)}Q}){_{(n)}\boldsymbol{e}}{_{(n)}\boldsymbol{e}_{j}}^{T},\ \ n\geq j.
\end{equation}
For this special augmentation, $_{(n)}E$ constitutes a close class for $_{(n)}\bm{\Phi}$.
Hence $_{(n)}\widetilde{Q}$ has a unique invariant probability vector.

In fact, $_{(n)}\widetilde{Q}$ in Remark \ref{the-difference} is also a linearly augmented truncation $q$-matrix  given by (\ref{aug-tru-2}) with $j=0$.
Thus, we cannot expect  the increasing property of the sequences $\{_{(n)}\delta_{j}\}$ or $\{_{(n)}\xi_{j}(|{_{(n)}}\bm{g}|)\}$
when $i=j$ for the linearly augmented truncation.

\begin{theorem}\label{the-cor-2}
Let $_{(n)}\widetilde{Q}$ be defined by (\ref{aug-tru-2}),
then both (i) and (ii) of Theorem \ref{the-cor-3-1} hold.
\end{theorem}

\proof  The first assertion is similar to that for Lemma \ref{lem-2} when $i\neq j$.
However, since the property in (\ref{sencor-pi-1}) does not hold for the linearly augmented truncation,
we cannot present a unified proof for both the censored Markov chain and the linearly augmented truncation.

Now, we prove the second assertion. It follows from $\bm{\pi}^{T}|\bm{g}|<\infty$ and Proposition \ref{tau-con-1} (ii) that
 $\mathbb{E}_{j}\left[{\xi_{j}(\bm{|g|})}\right]<\infty$.
By the strong Markov property, we have
\begin{eqnarray*}
   \mathbb{E}_{j}\left[{\xi_{j}(\bm{|g|})}\right]&=&
   \mathbb{E}_{j}\bigg[\int_{0}^{J_1}|{g}({\Phi}_t)|dt+\int_{J_1}^{\delta_{j}}|{g}({\Phi}_t)|dt\bigg] \\
   &=& \mathbb{E}_{j}\left[{J_1}|{g}(j)|\right]+\mathbb{E}_{j}\bigg[\int_{J_1}^{\delta_{j}}|{g}({\Phi}_t)|dt\bigg]\\
   &=& \frac{|{g}(j)|}{q_j}+\frac{1}{q_j}\sum_{k\in E,k\neq j}q_{jk}\mathbb{E}_{k}\left[{\xi_{j}(\bm{|g|})}\right],\\
\end{eqnarray*}
from which,
\begin{equation}\label{the-lin-poi-1}
  \sum_{k\in E,k\neq j}{q_{jk}}\mathbb{E}_{k}\left[{\xi_{j}(\bm{|g|})}\right]<\infty.
\end{equation}

Similarly, it follows that for any $n\geq j$,
\begin{equation}\label{the-lin-poi-2}
 \mathbb{E}_{j}\left[{_{(n)}\xi_{j}({_{(n)}}\bm{g})}\right]=\frac{{g}(j)}{{_{(n)}}{{\widetilde{q}}}_j}
 +\frac{1}{{_{(n)}}{{\widetilde{q}}}_j}\sum_{0\leq k\neq j\leq n}
   {q}_{jk}\mathbb{E}_{k}\left[{_{(n)}\xi_{j}({_{(n)}}\bm{g})}\right].
\end{equation}
Note that ${_{(n)}}{{\widetilde{q}}}_{jk}=q_{jk}$ for $0\leq k\neq j\leq n$.
By (\ref{the-lin-poi-1})--(\ref{the-lin-poi-2}) and the first assertion,
\begin{eqnarray}\label{the-lin-poi-3}
\nonumber&&\lim_{n\rightarrow\infty}\mathbb{E}_{j}\left[{_{(n)}\xi_{j}({_{(n)}}\bm{g})}\right]\\
\nonumber  &=& \lim_{n\rightarrow\infty}\frac{{g}(j)}{{_{(n)}}{{\widetilde{q}}}_j}+ \lim_{n\rightarrow \infty}\frac{1}{{_{(n)}}{{\widetilde{q}}}_j}\sum_{k\in E,k\neq j}
q_{jk}\left(\mathbb{E}_{k}\left[{_{(n)}\xi_{j}({_{(n)}}\bm{g}^{+})}\right]
-\mathbb{E}_{k}\left[{_{(n)}\xi_{j}({_{(n)}}\bm{g}^{-})}\right]\right)\mathbb{I}_{\{k\leq n\}}\\
\nonumber  &=&\frac{{g}(j)}{q_j}+ \frac{1}{q_j}\sum_{k\in E,k\neq j}
q_{jk}\left(\mathbb{E}_{k}\left[{\xi_{j}(\bm{g}^{+})}\right]
-\mathbb{E}_{k}\left[{\xi_{j}(\bm{g}^{-})}\right]\right)\\
 &=& \mathbb{E}_{j}\left[{\xi_{j}(\bm{g})}\right].
\end{eqnarray}
Similar to (\ref{the-lin-poi-3}),
\begin{equation}\label{the-lin-poi-5}
  \lim_{n\rightarrow\infty}\mathbb{E}_{j}\left[{{_{(n)}}{\delta}_{j}}\right]
  =\mathbb{E}_{j}\left[\delta_{j}\right].
\end{equation}
Thus, from Theorem 1.2 in Section 6.1 of \cite{a03}, and (\ref{the-lin-poi-3})--(\ref{the-lin-poi-5}), we then obtain
\begin{equation}\label{the-lin-poi-6}
\lim_{n\rightarrow\infty}{_{(n)}}\bm{\pi}^{T}{_{(n)}}\bm{g}
=\lim_{n\rightarrow\infty}\frac{\mathbb{E}_{j}\left[{_{(n)}\xi_{j}({_{(n)}}\bm{g})}\right]}
{\mathbb{E}_j[{_{(n)}}{\delta}_{j}]}
=\frac{\mathbb{E}_{j}\left[{\xi_{j}(\bm{g})}\right]}
{\mathbb{E}_j[\delta_{j}]}=\bm{\pi}^{T}\bm{g}.
\end{equation}
The rest of the proof follows from  Theorem \ref{the-cor-1}.
\BOX

Let $\bm{g}=\boldsymbol{e}_i$ in (\ref{the-lin-poi-6}),
we have the following interesting corollary directly,
which establishes the counterpart to the discrete-time result Theorem 3.2 in Diana and Seneta \cite{DS87}.

\begin{corollary}
Let $_{(n)}\widetilde{Q}$ be defined by (\ref{aug-tru-2}).
Then we have for any $i\in E$,
\[\lim_{n\rightarrow\infty}{_{(n)}}{\pi}(i)=\pi(i). \]
\end{corollary}

\begin{theorem}\label{variance-trunc-1}
Let $_{(n)}\widetilde{Q}$ be defined by (\ref{aug-tru-2}).
If $\mathbb{E}_\ell[\delta_\ell^2]<\infty$ and (\ref{var-cos-1}) holds for some $\ell\in E$, then we have
\[\lim_{n\rightarrow+\infty}{_{(n)}{\sigma}^{2}({_{(n)}}\bm{g})}=\sigma^{2}(\bm{g}).\]
\end{theorem}
\proof
Since we cannot perform the same arguments as that for (\ref{the-var-of-censor-1}),
we use the first expression of $\sigma^{2}(\bm{g)}$ in (\ref{var-exp}).
By the strong Markov property, it follows that
\begin{eqnarray*}
 \mathbb{E}_{j}\left[{\xi_{j}^{2}(\bm{|g|})}\right]
  &=& \mathbb{E}_{j}\bigg[\bigg(\int_0^{J_1}|g(\Phi_t)|dt+\int_{J_1}^{\delta_{j}}|g(\Phi_t)|dt\bigg)^2\bigg]\\
  &=& \mathbb{E}_{j}\bigg[\bigg({J_1}|g(j)|\bigg)^2\bigg]
  +2\mathbb{E}_{j}\bigg[{J_1}|g(i)|\bigg(\int_{J_1}^{\delta_{j}}|g(\Phi_t)|dt\bigg)\bigg]
  +\mathbb{E}_{j}\bigg[\bigg(\int_{J_1}^{\delta_{j}}|g(\Phi_t)|dt\bigg)^2\bigg]\\
  &=&\frac{2|g(j)|^{2}}{q_j^{2}}
  +\frac{2|g(j)|}{q_j^{2}}\sum_{k\neq j, k\in E}q_{jk}\mathbb{E}_{k}\left[{\xi_{j}(\bm{|g|})}\right]
  +\frac{1}{q_j}\sum_{k\neq j, k\in E}q_{jk}\mathbb{E}_{k}\left[{\xi_{j}^{2}(\bm{|g|})}\right],
\end{eqnarray*}
from which,
\begin{equation}\label{the-lin-var-1}
\sum_{k\neq j, k\in E}q_{jk}\mathbb{E}_{k}\left[{\xi_{j}^{2}(\bm{|g|})}\right]<\infty.
\end{equation}

Similar to DTMCs,
\begin{eqnarray*}
 \mathbb{E}_{j}\left[{_{(n)}\xi_{j}^{2}({_{(n)}}\bm{\overline{g}})}\right]
 = \mathbb{E}_{j}\left[{_{(n)}\xi_{j}^{2}({_{(n)}}\bm{g})}\right]
-2\left({_{(n)}}\bm{\pi}^{T}{_{(n)}}\bm{g}\right)
\mathbb{E}_{j}\left[{_{(n)}\xi_{j}({_{(n)}}\bm{g})}\cdot{_{(n)}}\delta_{j}\right]
+\left({_{(n)}}\bm{\pi}^{T}{_{(n)}}\bm{g}\right)^{2}\mathbb{E}_{j}[{_{(n)}}\delta_{j}^{2}].
\end{eqnarray*}
For any $n\geq j$, we obtain that
\begin{eqnarray}\label{the-lin-var-2}
\nonumber\mathbb{E}_{j}\left[{_{(n)}\xi_{j}^{2}({_{(n)}}\bm{g})}\right] &=&
\frac{2g(j)^{2}}{{_{(n)}}{\widetilde{q}_j}^{2}}
  +\frac{2g(j)}{{_{(n)}}{\widetilde{q}}_j^{2}}\sum_{0\leq k\neq j\leq n}
  {q}_{jk}\mathbb{E}_{k}\left[{_{(n)}\xi_{j}({_{(n)}}\bm{g})}\right]
+\frac{1}{{_{(n)}}{\widetilde{q}}_j}\sum_{0\leq k\neq j\leq n}q_{jk}
  \mathbb{E}_{k}\left[{_{(n)}\xi_{j}^{2}({_{(n)}}\bm{g})}\right].
\end{eqnarray}
It follows from (\ref{the-lin-var-1}) and Theorem \ref{the-cor-2} that
\begin{eqnarray}\label{the-lin-var-4}
\nonumber&&\lim_{n\rightarrow\infty}\sum_{0\leq k\neq j\leq n}q_{jk}
  \mathbb{E}_{k}\left[{_{(n)}\xi_{j}^{2}({_{(n)}}\bm{g})}\right]\\
\nonumber  &=&  \lim_{n\rightarrow \infty} \sum_{0\leq k\neq j\leq n}q_{jk}\mathbb{E}_{k}\left[{_{(n)}\xi_{j}^{2}({_{(n)}}\bm{g}^{+})}\right]
 +\lim_{n\rightarrow \infty} \sum_{0\leq k\neq j\leq n}
 q_{jk}\mathbb{E}_{k}\left[{_{(n)}\xi_{j}^{2}({_{(n)}}\bm{g}^{-})}\right]\\
\nonumber  &&-2\lim_{n\rightarrow \infty} \sum_{0\leq k\neq j\leq n}
q_{jk}\mathbb{E}_{k}\left[{_{(n)}\xi_{j}({_{(n)}}\bm{g}^{+})}{_{(n)}\xi_{j}({_{(n)}}\bm{g}^{-})}\right]\\
 &=& \sum_{k\neq j, k\in E}q_{jk}\mathbb{E}_{k}\left[{\xi_{j}^{2}(\bm{g})}\right].
\end{eqnarray}
Thus, by (\ref{the-lin-poi-3}) and (\ref{the-lin-var-1})--(\ref{the-lin-var-4}),
\begin{equation}\label{the-lin-var-5}
  \lim_{n\rightarrow\infty}\mathbb{E}_{j}\left[{_{(n)}\xi_{j}^{2}({_{(n)}}\bm{g})}\right]
  =\mathbb{E}_{j}\left[{\xi_{j}^{2}(\bm{g})}\right].
\end{equation}
Similarly, we have
\begin{equation}\label{the-lin-var-6}
 \lim_{n\rightarrow \infty} \mathbb{E}_{j}[{_{(n)}}{\delta}_{j}^{2}] =\mathbb{E}_{j}[\delta_{j}^{2}],
\end{equation}
and
\begin{equation}\label{the-lin-var-7}
    \lim_{n\rightarrow \infty} \mathbb{E}_{j}\left[{_{(n)}\xi_{j}^{2}({_{(n)}}\bm{g})}\cdot{_{(n)}}{\delta}_{j}\right] = \mathbb{E}_{j}\left[{\xi_{j}^{2}(\bm{g})}\cdot{\delta}_{j}\right].
\end{equation}
Hence, by the above  limits (\ref{the-lin-var-5})--(\ref{the-lin-var-7}) and (\ref{the-lin-poi-5}),
we obtain the assertion of this theorem.
\BOX

\section{Applications}

In this section, we apply our results to two types of classical asymmetric Markov processes:
single-birth processes and single-death processes.
For simplicity of the presentation, we only consider the discrete-time case for the former and the continuous-time case for the latter.

\subsection{Discrete-time single-birth processes}
The discrete-time single-birth process $\boldsymbol{\Phi}=\{\Phi_k: k\in \mathbb{Z}_+\}$
 has a special transition matrix given by $p_{i,i+1}>0$
and $p_{i,i+k}=0$ for all $i\geq0$ and $k\geq2$.
In order to solve Poisson's equation, we need the following notations (see, e.g. Chen and Zhang \cite{cz04}). Define
\[
p_{m}^{(k)}=\sum_{i=0}^{k}p_{mi},\ \ 0\leq k<m,
\]
and
\[
F_{i}^{(i)}=1, \ \ F_{m}^{(i)}=\frac{1}{p_{m,m+1}}\sum_{k=i}^{m-1}p_{m}^{(k)}F_{k}^{(i)}, \ 0\leq i<m.
\]

\begin{theorem}\label{infinte-poi-equ-0}
Suppose that the single-birth process $\boldsymbol{\Phi}$ is irreducible and positive recurrent.
If the function $\bm{g}$ satisfies $\bm{\pi}^{T}|\bm{g}|<\infty$, then for any fixed state $j\in E$, we have
\begin{equation}\label{single-birth-theorem}
  f_{j}(i)=
  \left\{\aligned
  & \sum_{m=i}^{j-1}\sum_{k=0}^{m}\frac{F_{m}^{(k)}\overline{g}(k)}{p_{k,k+1}}, &  if\  i<j,\\
  & 0, & if\  i= j,\\
  & -\sum_{m=j}^{i-1}\sum_{k=0}^{m}\frac{F_{m}^{(k)}\overline{g}(k)}{p_{k,k+1}}, & if\  j<i.
  \endaligned \right.
\end{equation}
\end{theorem}
\proof For a fixed state $j\in E$ and $n> j$,
we consider the $(j+1)$th augmented matrix $_{(n)}\widetilde{P}$ defined by (\ref{gen-tru-2-1}).
In the follows we solve (\ref{poi-fun-2}) with ${_{(n)}}f_j(j)=0$.
From Corollary 2.3 in \cite{cz04}, we have that
\[{_{(n)}}f_j(i)-{_{(n)}}f_j(i+1)=\sum_{k=0}^{i}\frac{F_{i}^{(k)}{_{(n)}}\overline{g}(k)}{p_{k,k+1}},\ \ 0\leq i\leq n-1.\]
Since $_{(n)}f_{j}(j)=0$,  we can use the inductive arguments to show for $i<j$,
\begin{equation}\label{single-birth-tru-1}
{_{(n)}}f_j(i)=\sum_{m=i}^{j-1}\sum_{k=0}^{m}\frac{F_{m}^{(k)}{_{(n)}}\overline{g}(k)}{p_{k,k+1}}.
\end{equation}
For the case of $j<i\leq n$,
\begin{equation}\label{single-birth-tru-2}
{_{(n)}}f_j(i)=-\sum_{m=0}^{i-1}\sum_{k=0}^{m}\frac{F_{m}^{(k)}{_{(n)}}\overline{g}(k)}{p_{k,k+1}}.
\end{equation}
Thus we obtain the solution of Poisson's function (\ref{poi-fun-2}).
It follows from Theorem \ref{the-cor-1} that for $i<j$,
\begin{eqnarray*}
  \lim_{n\rightarrow\infty}{_{(n)}}f_{j}(i)
   &=& \lim_{n\rightarrow\infty}\sum_{m=0}^{i-1}\bigg(\sum_{k=0}^{m}\frac{F_{m}^{(k)}g(k)}{p_{k,k+1}}
   -\left({_{(n)}}\bm{\pi}^{T}{_{(n)}}\bm{g}\right)\sum_{k=0}^{m}\frac{F_{m}^{(k)}}{p_{k,k+1}}\bigg)\\
   &=& \sum_{m=0}^{i-1}\sum_{k=0}^{m}\frac{F_{m}^{(k)}g(k)}{p_{k,k+1}}
   -\lim_{n\rightarrow\infty}\left({_{(n)}}\bm{\pi}^{T}{_{(n)}}\bm{g}\right)\sum_{k=0}^{m}\frac{F_{m}^{(k)}}{p_{k,k+1}}\\
   &=& \sum_{m=0}^{i-1}\sum_{k=0}^{m}\frac{F_{m}^{(k)}g(k)}{p_{k,k+1}}
   -\left(\bm{\pi}^{T}\bm{g}\right)\sum_{k=0}^{m}\frac{F_{m}^{(k)}}{p_{k,k+1}}\\
  &=& \sum_{m=0}^{i-1}\sum_{k=0}^{m}\frac{F_{m}^{(k)}\overline{g}(k)}{p_{k,k+1}}.
\end{eqnarray*}
The case of $i>j$ can be verified similarly, which is omitted here.
We obtain the assertion.
\BOX

\begin{remark}
The result of Theorem \ref{infinte-poi-equ-0} was first established in \cite{jly14}.
Here, we revisit it using the  truncation approximation.
The continuous-time case, which was first presented in \cite{lz15}, can be investigated similarly.
\end{remark}

\subsection{Continuous-time single-death processes}

We call $\mathbf{\Phi}=\{\Phi_t: t\in \mathbb{R}_+\}$ a single-death process on the state space $E$,
if its $q$-matrix $Q$ satisfies $q_{i,i-1}>0$ for all $i\geq1$ and $q_{i,i-k}=0$ for all $i\geq k\geq2$.
Assume that $Q$ is a totally stable  and regular single-death $q$-matrix.
To solve Poisson's equation, we need to introduce the following notations (see e.g. Zhang \cite{z18}):
\[
q_{m}^{(k)}=\sum_{i=k}^{\infty}q_{mi},\ \ k>m\geq0,
\]
and
\[
G_{i}^{(i)}=1, \ \ G_{m}^{(i)}=\frac{1}{q_{m,m-1}}\sum_{k=m+1}^{i}q_{m}^{(k)}G_{k}^{(i)}, \ 1\leq m<i.
\]

It is shown in Liu et al. \cite{lt15} that the $q$-matrix $_{(n)}\widetilde{Q}$ of the censored Markov chain on $_{(n)}E$
is actually the last column augmentation of $Q$.
Specifically,  $_{(n)}\widetilde{Q}$ is given by
\[
_{(n)}\widetilde{Q}=\left (
 \begin{array}{cccccc}
  q_{00} & q_{01} & q_{02} &\cdots & q_{0,n-1}& q_{0}^{(n)} \\
  q_{10} & q_{11} & q_{12} &\cdots & q_{1,n-1}& q_{1}^{(n)}\\
  0&q_{21} & q_{22} & \cdots & q_{2,n-1}& q_{2}^{(n)}\\
  \vdots&\vdots&\vdots&\ddots& \vdots& \vdots \\
  0& 0& 0& \cdots&q_{n-1,n-1}&q_{n-1}^{(n)} \\
   0&  0&0& \cdots&q_{n,n-1}& -q_{n,n-1}\
   \end{array}\right ).
 \]

\begin{theorem}\label{infinte-poi-equ}
Suppose that the single-death process $\mathbf{\Phi}$ is irreducible and positive recurrent.
If the function $\bm{g}$ satisfies $\bm{\pi}^{T}|\bm{g}|<\infty$,
then for any fixed state $j\in E$, we have
\begin{equation}\label{sin-dea-equ}
 f_{j}(i)=\left\{\aligned & -\sum_{m=i+1}^{j}\sum_{k=m}^{\infty}\frac{G_{m}^{(k)}\overline{g}(k)}{q_{k,k-1}}, &   if\  i<j,\\
  &0, & if\  i=j,\\
  &\sum_{m=j+1}^{i}\sum_{k=m}^{\infty}\frac{G_{m}^{(k)}\overline{g}(k)}{q_{k,k-1}}, & if\  j<i. \endaligned \right.
\end{equation}
Moreover, if (\ref{var-cos-1}) holds for some $\ell$, we have
\begin{equation}\label{the-var-cos}
  \sigma^{2}(\bm{g})=2\sum_{i=1}^{\infty}\pi(i) \bar{g}(i)\sum_{m=1}^{i}\sum_{k=m}^{\infty}\frac{G_{m}^{(k)}\overline{g}(k)}{q_{k,k-1}}.
\end{equation}

\end{theorem}
\proof
Similar to the analysis of the single-birth processes,
from Poisson's equation (\ref{poi-fun-3}) and for a fixed state $j<n$, we obtain
\begin{equation*}\label{p-1-2}
  {_{(n)}}f_{j}(i)-{_{(n)}}f_{j}(i-1)=\frac{1}{q_{i,i-1}}\bigg(\sum_{k=i+1}^{n}{q}_{i}^{(k)}
  \big({_{(n)}}f_{j}(k)-{_{(n)}}f_{j}(k-1)\big)+{_{(n)}}\bar{g}(i)\bigg), \ \ 1\leq i\leq n-1.
\end{equation*}
\begin{equation*}\label{p-1-3}
  {_{(n)}}f(n)-{_{(n)}}f(n-1)=\frac{{_{(n)}}\overline{g}(n)}{q_{n,n-1}}.
\end{equation*}
From Corollary 2.3 in \cite{z18}, one can easily show that
\begin{equation}\label{p-1-4}
 {_{(n)}}f_{j}(i)-{_{(n)}}f_{j}(i-1)=\sum_{k=i}^{n}\frac{G_{i}^{(k)}{_{(n)}}\overline{g}(k)}{q_{k,k-1}},\ \ 1\leq i\leq n.
\end{equation}
Since $_{(n)}f(j)=0$, by using the induction similar to (\ref{single-birth-tru-1}) and (\ref{single-birth-tru-2}),  it follows that
\[{_{(n)}}f(i)=-\sum_{m=i+1}^{j}\sum_{k=m}^{n}\frac{G_{m}^{(k)}{_{(n)}}\overline{g}(k)}{q_{k,k-1}},  \ \ 0\leq i<j,\]
and
\[
{_{(n)}}f(i)=\sum_{m=j+1}^{i}\sum_{k=m}^{n}\frac{G_{m}^{(k)}{_{(n)}}\overline{g}(k)}{q_{k,k-1}},\ \  j<i\leq n.
\]
Then, we obtain the solution of Poisson's function (\ref{poi-fun-3}).

From Theorem \ref{the-cor-3-1}, we have, for $i<j$
\begin{eqnarray*}
\lim_{n\rightarrow\infty}{_{(n)}f_{j}(i)}
&=&-\sum_{m=i+1}^{j}\lim_{n\rightarrow\infty}\sum_{k=m}^{n}\frac{{G}_{m}^{(k)}{_{(n)}\overline{g}(k)}}{{q}_{k,k-1}}\\
&=&-\sum_{m=i+1}^{j}\bigg(\lim_{n\rightarrow\infty}\sum_{k=m}^{n}\frac{G_{m}^{(k)}g(k)}{q_{k,k-1}}
-\lim_{n\rightarrow\infty}\left({_{(n)}\bm{\pi}^{T}{_{(n)}}\bm{g}}\right)\sum_{k=m}^{n}\frac{G_{m}^{(k)}}{q_{k,k-1}}\bigg)\\
&=&-\sum_{m=i+1}^{j}\bigg(\sum_{k=m}^{\infty}\frac{G_{m}^{(k)}g(k)}{q_{k,k-1}}
-\left(\bm{\pi}^{T}\bm{g}\right)\sum_{k=m}^{\infty}\frac{G_{m}^{(k)}}{q_{k,k-1}}\bigg)\\
&=&-\sum_{m=i+1}^{j}\sum_{k=m}^{\infty}\frac{G_{m}^{(k)}\overline{g}(k)}{q_{k,k-1}}.
\end{eqnarray*}
The case of $i>j$ can be verified similarly. Thus, the first assertion is proved.

According to (\ref{var-exp}) and (\ref{sin-dea-equ}), we have for any $j\in E$,
\begin{eqnarray*}
\sigma^{2}(\bm{g})
&=&-2\sum_{i=0}^{j-1}\pi(i)\overline{g}(i)\sum_{m=i+1}^{j}\sum_{k=m}^{\infty}\frac{G_{m}^{(k)}\overline{g}(k)}{q_{k,k-1}}
+2\sum_{i=j+1}^{\infty}\pi(i) \overline{g}(i)\sum_{m=j+1}^{i}\sum_{k=m}^{\infty}\frac{G_{m}^{(k)}\overline{g}(k)}{q_{k,k-1}} \\
&=&2\sum_{i=1}^{\infty}\pi(i) \overline{g}(i)\sum_{m=1}^{j}\sum_{k=m}^{\infty}\frac{G_{m}^{(k)}\overline{g}(k)}{q_{k,k-1}}
-2\sum_{i=1}^{j-1}\pi(i) \overline{g}(i)\sum_{m=i+1}^{j}\sum_{k=m}^{\infty}\frac{G_{m}^{(k)}\overline{g}(k)}{q_{k,k-1}}\\
&&+2\sum_{i=j+1}^{\infty}\pi(i)\overline{g}(i)\sum_{m=j+1}^{i}\sum_{k=m}^{\infty}\frac{G_{m}^{(k)}\overline{g}(k)}{q_{k,k-1}}\\
&=&2\sum_{i=1}^{\infty}\pi(i) \overline{g}(i)\sum_{m=1}^{i}\sum_{k=m}^{\infty}\frac{G_{m}^{(k)}\overline{g}(k)}{q_{k,k-1}},
\end{eqnarray*}
in which the third equation following using the fact that
\[\sum_{i=0}^{\infty}\pi(i)\overline{g}(i)=0.\]
\BOX

\begin{remark} \
(i) The integral-type functionals for single death processes had been investigated by Wang and Zhang \cite{WZ20}
by using different arguments, which hold only for downward integral-type functionals (i.e. the case $i>j$ in (\ref{sin-dea-equ})).

(ii) When $\mathbf{\Phi}$ is an irreducible and positive recurrent birth-death process, (\ref{the-var-cos}) becomes
\[
  \sigma^{2}(\bm{g})=\sum_{i=0}^{\infty}\frac{1}{q_{i,i+1}\pi(i)}\bigg(\sum_{k=0}^{i}\pi(k)\overline{g}(k)\bigg)^{2}.
\]
\end{remark}

\begin{example}\label{exm-2}
This example was taken from \cite{z18}.
Give a constant $b>2$. Define a totally stable, conservative,
and irreducible single-death $q$-matrix $Q=(q_{ij})$ as follows:
\[q_{ij}=\frac{b-1}{b^{j-i+1}},\ \ j\geq i+1;\ \ q_{i,i-1}=\frac{b-1}{b},\ \ q_{i}=-q_{ii}=\frac{b^{2}-b+1}{b^{2}},\ \ i\geq1;\]
\[q_{0j}=\frac{b-1}{b^{j+1}},\ \ j\geq 1;\ \ q_0=-q_{00}=\frac{1}{b}.\]
\end{example}
By calculations, we have
\[q_{n}^{(k)}=\frac{1}{b^{k-n+1}},\ \ 1\leq n\leq k;\ \ \ q_{0}^{(k)}=\frac{1}{b^{k}},\ \ k\geq1;\]
and
\[G_{n}^{(i)}=\frac{1}{b(b-1)^{i-n}},\ \ 1\leq n\leq i.\]

From \cite{z18} again, we know that the single-death process is exponentially ergodic,
and the unique invariant probability vector is given by
\[\pi(i)=\frac{b-2}{(b-1)^{i+1}},\ \ i\geq0.\]
Let $g(i)=i$, then
\[\bm{\pi}^{T}|\bm{g}|=\bm{\pi}^{T}\bm{g}=\frac{1}{b-2},\ \ \bm{\pi}^{T}|\bm{g}|^{3}=\frac{b^{2}+2b-2}{(b-1)^{3}}.\]
Thus, from Theorem \ref{infinte-poi-equ}, we have
\begin{equation*}\label{exm-poi-sol}
  f_{j}(i)=\left\{\aligned &\frac{(i-j)[(i+j+1)(b-1)-2]}{2(b-2)},& if\  i\neq j,\\
  &0, & if\  i= j. \endaligned \right.
\end{equation*}
According to Remark \ref{the-condition-constant}, we know that
the variance constant $\sigma^{2}(\bm{g})$ exists and is given by
\begin{equation*}
  \sigma^{2}(\bm{g})=\frac{2b^{3}-6b^{2}+8b-4}{(b-2)^{4}}.
\end{equation*}

\begin{example}\label{exm-3} Consider an extended class of branching processes with $q$-matrix as follows:
\begin{equation*}\label{branching-model}
  q_{ij}=\left\{\aligned
  & b_j, &i=0,j\geq1,\\
  &i^{\alpha}p_{j-i+1}, &i\geq 1,\ j\geq i-1, \\
  &0,& else,
  \endaligned \right.
\end{equation*}
where $\alpha>0$, $b_j\geq0$ for $j\geq1$ and $0<-b_0=\sum_{k\neq0}b_j<\infty$; $p_j\geq0$ for $j\neq1$, and $0<-p_1=\sum_{k\neq0}p_j<\infty$.
\end{example}
In the above model, let
\[\alpha=1,\ \ b_{j}=\left(\frac{1}{3}\right)^{j-1}\frac{2}{3},\ \ j\geq1;\ \ p_{j}=\left(\frac{1}{3}\right)^{j},\ \ j\neq1.\]
By calculations, we know
\[q_{n}^{(k)}=\frac{n}{2\cdot3^{k-n}},\ \ 1\leq n\leq k;\ \ \ q_{0}^{(k)}=\frac{1}{3^{k-1}},\ \ k\geq1;\]
and
\[G_{n}^{(i)}=\frac{1}{3\cdot2^{i-n}},\ \ 1\leq n\leq i.\]
Let
\[h_{n}=\sum_{k=n}^{\infty}\frac{G_{n}^{(k)}}{q_{k,k-1}}=\frac{1}{n}+\frac{2^{n}}{3}\left(\ln{2}-\sum_{k=1}^{n}\frac{1}{k2^{k}}\right),\ \ n\geq1.\]

From Theorem 4.1 in \cite{z18}, we know that the single-death process is exponentially ergodic.
According to \cite{lt15}, the unique invariant probability vector is given by
\[\pi(0)=\frac{1}{1+\ln{4}},\ \ \pi(i)=\frac{1}{i\cdot2^{i-1}(1+\ln{4})},\ \ i\geq1.\]
Let $g(i)=i$, then
\[\bm{\pi}^{T}|\bm{g}|=\bm{\pi}^{T}\bm{g}=\frac{2}{1+\ln{4}},\ \ \bm{\pi}^{T}|\bm{g}|^{3}=\frac{12}{1+\ln{4}}.\]
Thus, from Theorem \ref{infinte-poi-equ}, we obtain
\begin{equation*}\label{exm-poi-sol-2}
  f_{j}(i)=\left\{\aligned &\frac{4}{3}(i-j)+\frac{2}{1+\ln{4}}\sum_{n=i+1}^{j}h_n,& if\  i<j,\\
  &0, & if\  i= j. \\
  & \frac{4}{3}(i-j)-\frac{2}{1+\ln{4}}\sum_{n=j+1}^{i}h_n,&if\  i> j.\endaligned \right.
\end{equation*}
According to Remark \ref{the-condition-constant}, 
we know that the variance constant $\sigma^{2}(\bm{g})$ exists and is given by
\begin{equation}\label{example-5-2}
\sigma^{2}(\bm{g})=\frac{2}{1+\ln{4}}\sum_{i=1}^{\infty}\frac{1}{i2^{i-1}}\left(i-\frac{2}{1+\ln{4}}\right)
  \left(\frac{4}{3}i-\frac{2}{1+\ln{4}}\sum_{n=1}^{i}h_n\right).
\end{equation}

On the one hand, we can approximate $\sigma^{2}(\bm{g})$
by truncating the corresponding infinite series (\ref{example-5-2}). Denote by $_{(n)}\sigma^{2}(\bm{g})$
the partial sum of the first $n$ items in the series (\ref{example-5-2}).
To gain information about $n$, we bound the error between $\sigma^{2}(\bm{g})$ and $_{(n)}\sigma^{2}(\bm{g})$ as follows
\[
  e_n:=\sigma^{2}(\bm{g})-{_{(n)}\sigma^{2}(\bm{g})} \leq \sum_{i=n}^{\infty}\frac{i}{2^{i-2}}.
\]
Moreover, $e_n\leq 10^{-4}$ when $n>22$. On the other hand,
according to Theorem \ref{variance-trunc-2}, we can approximate $\sigma^{2}(\bm{g})$ with $_{(n)}\sigma^{2}({_{(n)}\bm{g}})$ directly.
The comparison between $_{(n)}\sigma^{2}(\bm{g})$ and $_{(n)}\sigma^{2}({_{(n)}\bm{g}})$ is depicted in Table 1,
which shows both are almost identical and the variance constant is 1.4645.

\begin{table}[h]
\caption{The variation of $_{(n)}\sigma^{2}(\bm{g})$ and $_{(n)}\sigma^{2}({_{(n)}\bm{g}})$ with the level $n$.}
  \centering
  \begin{tabular}{ccccccccccccc}
    \hline
    $n$ & 10 & 12&14&16&18&20&22&24&26\\
    \hline
$_{(n)}\sigma^{2}(\bm{g})$          & 1.4448 & 1.4585& 1.4627 & 1.4640 & 1.4643 & 1.4644 & 1.4645 & 1.4645 & 1.4645\\
$_{(n)}\sigma^{2}({_{(n)}\bm{g}})$  & 1.4394 & 1.4566& 1.4621 & 1.4638 & 1.4643 & 1.4644 & 1.4644 & 1.4645 & 1.4645\\
    \hline
  \end{tabular}
\end{table}

\section{Conclusion and discussion}

We develop the technique of augmented truncation approximations for the solution of Poisson's equation and the variance constants in CLTs. The role of the technique is two-fold. On the one hand,  it provides us a useful way to
derive explicit expressions for the solution and the variance constant for some infinite-state Markov chains.
On the other hand, it provide us an efficient way to approximate them  numerically as the truncation size becomes large.

We now discuss possible extensions and improvements of the results in this paper.

It is interesting to extend  the technique of augmented truncation  to investigate block-structured Markov chains. The censored chain technique can be expected to hold with a little  more complicated arguments. The extension of the technique of linearly augmented truncation is more involved since first it should be extended to block column augmentation, and meanwhile, the monotone assumption about the first return time moments should be extended to block monotone situation.

To perform the truncation approximations effectively, it is desirable to know some information about the truncation size. Hence it is important to investigate the bounds on the  truncation approximation error. The arguments in this paper and some ideas in Liu and Li \cite{ll18} may be used, but definitely it requires also new  arguments. This is an interesting and challenging topic for the future research.

\appendix
\section{ }
We present two useful propositions about the first return time for DTMCs and CTMCs, respectively.

\begin{proposition}\label{tau-con-0}
Suppose that the DTMC $\mathbf{\Phi}$ is irreducible and positive recurrent.
For any finite non-negative vector $\bm{g}$, we have
\begin{description}
  \item[(i)] $\mathbb{E}_{\ell}[\zeta_{\ell}^{p}(\bm{g})]<\infty$ for some $\ell\in E$
if and only if $\mathbb{E}_{i}[\zeta_{i}^{p}(\bm{g})]<\infty$ for any state $i\in E$;
  \item[(ii)] $\mathbb{E}_{\ell}[\zeta_{\ell}(\bm{g})]<\infty$ for some $\ell\in E$
if and only if $\mathbb{E}_{i}[\zeta_{j}(\bm{g})]<\infty$ for any $i,j\in E$.
In particular, if $\bm{\pi}^{T}\bm{g}<\infty$,
then $\mathbb{E}_{i}[\zeta_{j}(\bm{g})]<\infty$ for any $i,j\in E$.
\end{description}
\end{proposition}

\proof By Theorem 4 of Section 14 in Chung \cite{c60}, we obtain the first assertion directly.
For the second assertion,
we only need to prove the sufficiency since the necessity is obvious.
Since the chain is irreducible, then there exists some $m>0$ such that
$$
_{i}p_{ij}^{m}:=P\{\Phi_m=j, \Phi_k\neq i, 1\leq k\leq m|\Phi_0=i\}>0.
$$
Thus, we have
\begin{eqnarray*}
  \mathbb{E}_{i}\bigg[\sum_{k=0}^{\tau_{i}-1}{g}(\Phi_k)\bigg]
   &=&\mathbb{E}_{i}\bigg[ \sum_{k=1}^{\infty}g(\Phi_k)\mathbb{I}_{\{k<\tau_{i}\}}\bigg]+g(i)\\
   &=&\sum_{\ell\in E}\sum_{k=1}^{\infty}{_{i}p_{i\ell}^{k}}g(\ell)+g(i)\\
   &\geq&\sum_{\ell\in E}\sum_{k=m+1}^{\infty}\sum_{j\in E}{_{i}p_{ij}^{m}}\cdot{_{i}p_{j\ell}^{k-m}}g(\ell)\\
   &\geq&{_{i}p_{ij}^{m}}\sum_{\ell\in E}\sum_{k=m+1}^{\infty} {_{i}p_{j\ell}^{k-m}}g(\ell)\\
   &=&{_{i}p_{ij}^{m}}\sum_{\ell\in E}\sum_{k=1}^{\infty} {_{i}p_{j\ell}^{k}}g(\ell),
\end{eqnarray*}
from which,
\[\mathbb{E}_{j}[\zeta_{i}(\bm{g})]=\mathbb{E}_{j}\bigg[\sum_{k=0}^{\tau_{i}-1}{g}(\Phi_k)\bigg]
=\sum_{\ell\in E}\sum_{k=1}^{\infty} {_{i}p_{j\ell}^{k}}g(\ell)+g(j)<\infty.\]
Since $\mathbf{\Phi}$ is positive recurrent, this shows that $\mathbb{E}_{l}[\delta_{l}]<\infty$ for any $\ell\in E$. If $\bm{\pi}^{T}\bm{g}<
\infty$, then by Theorem 10.0.1 in \cite{mt09}, we obtain
\begin{equation*}\label{sta-funct}
\mathbb{E}_{l}[\zeta_{l}(\bm{g})]= \mathbb{E}_{l}\bigg[\sum_{k=0}^{\tau_{l}-1}{g}(\Phi_k)\bigg]
= \left(\bm{\pi}^{T}\bm{g}\right) \mathbb{E}_{l}[\delta_{l}]<\infty.
\end{equation*}
Thus the proof is finished.
\BOX

\begin{proposition}\label{tau-con-1}
Suppose that the CTMC $\mathbf{\Phi}$ is irreducible and positive recurrent.
For any non-negative finite function $\bm{g}$, we have
\begin{description}
  \item[(i)] $\mathbb{E}_{\ell}[\xi_{\ell}^{p}(\bm{g})]<\infty$ for some $\ell\in E$
if and only if $\mathbb{E}_{i}[\xi_{i}^{p}(\bm{g})]<\infty$ for any state $i\in E$;
  \item[(ii)] $\mathbb{E}_{\ell}[\xi_{\ell}(\bm{g})]<\infty$ for some $\ell\in E$
if and only if $\mathbb{E}_{i}[\xi_{j}(\bm{g})]<\infty$ for any $i,j\in E$.
In particular, if $\bm{\pi}^{T}\bm{g}<\infty$,
then $\mathbb{E}_{i}[\zeta_{j}(\bm{g})]<\infty$ for any $i,j\in E$.
\end{description}
\end{proposition}
\proof The proof is similar to Proposition \ref{tau-con-0}, which is omitted here.
\BOX

\section*{Acknowledgements}
This research was supported in part by 
the National Natural Science Foundation of China (Grants No. 11971486, 11771452),
Natural Science Foundation of Hunan (Grants No. 2019JJ40357, 2020JJ4674), 
the Innovation Program of Central South University (Grant No. 2020zzts039), and
the Natural Sciences and Engineering Research Council (NSERC) of Canada (Discovery Grant).




\begin{thebibliography}{99}

\bibitem{mt09}
S.~Meyn, R.~Tweedie, Markov Chains and Stochastic Stability, 2nd Edition,
  Cambridge: Cambridge University Press, 2009.
\newblock \href {http://dx.doi.org/10.1007/978-1-4471-3267-7}
  {\path{doi:10.1007/978-1-4471-3267-7}}.

\bibitem{b07}
D.~Bertsekas, Dynamic Programming and Optimal Control, 3rd Edition, Cambridge:
  Athena Scientific, 2007.

\bibitem{gm96}
P.~Glynn, S.~Meyn, A liapounov bound for solutions of the poisson equation,
  Ann. Probab. 24 (1996) 916--931.
\newblock \href {http://dx.doi.org/10.1214/aop/1039639370}
  {\path{doi:10.1214/aop/1039639370}}.

\bibitem{l15}
Y.~Liu, Perturbation analysis for continuous-time markov chains, Sci.
  China-Math. 58(12) (2015) 2633--2642.
\newblock \href {http://dx.doi.org/10.1007/s11425-015-5019-z}
  {\path{doi:10.1007/s11425-015-5019-z}}.

\bibitem{sm94}
S.~Asmussen, M.~Bladt, Poisson's equation for queues driven by a markovian
  marked point process, Queueing Syst. 17 (1994) 235--274.
\newblock \href {http://dx.doi.org/10.1007/BF01158696}
  {\path{doi:10.1007/BF01158696}}.

\bibitem{go02}
P.~Glynn, P.~Ormoneit, Hoeffding's inequality for uniformly ergodic markov
  chains, Statist. Probab. Lett. 56 (2002) 143--146.
\newblock \href {http://dx.doi.org/10.1016/S0167-7152(01)00158-4}
  {\path{doi:10.1016/S0167-7152(01)00158-4}}.

\bibitem{CL19}
M.~Choi, E.~Li, A hoeffding's inequality for uniformly ergodic diffusion
  process, Statist. Probab. Lett. 150 (2019) 23--28.
\newblock \href {http://dx.doi.org/10.1016/j.spl.2019.02.012}
  {\path{doi:10.1016/j.spl.2019.02.012}}.

\bibitem{ll18}
Y.~Liu, W.~Li, Error bounds for augmented truncation approximations of markov
  chains via the perturbation method, Adv. Appl. Probab. 50(2) (2018) 645--669.
\newblock \href {http://dx.doi.org/10.1017/apr.2018.28}
  {\path{doi:10.1017/apr.2018.28}}.

\bibitem{MS02}
A.~Makowski, A.~Shwartz, The Poisson equation for countable Markov chains:
  probabilistic methods and interpretations, Handbook of Markov Decision
  Processes. Springer US, 2002.
\newblock \href {http://dx.doi.org/10.1007/978-1-4615-0805-2_9}
  {\path{doi:10.1007/978-1-4615-0805-2_9}}.

\bibitem{m15}
H.~Masuyama, Error bounds for augmented truncations of discrete-time
  block-monotone markov chains under geometric drift conditions, Adv. Appl.
  Probab. 47 (2015) 83--105.
\newblock \href {http://dx.doi.org/10.1239/aap/1427814582}
  {\path{doi:10.1239/aap/1427814582}}.

\bibitem{m16}
H.~Masuyama, Error bounds for augmented truncations of discrete-time
  blockmonotone markov chains under subgeometric drift conditions, SIAM J.
  Martix Anal. Appl 37 (2016) 877--910.
\newblock \href {http://dx.doi.org/10.1137/15M1024743}
  {\path{doi:10.1137/15M1024743}}.

\bibitem{l10}
Y.~Liu, Augmented truncation approximations of discrete-time markov chains,
  Oper. Res. Lett. 38(3) (2010) 218--222.
\newblock \href {http://dx.doi.org/10.1016/j.orl.2009.12.001}
  {\path{doi:10.1016/j.orl.2009.12.001}}.

\bibitem{J04}
Jones, L.~Galin, On the markov chain central limit theorem, Probab. Surv. 1
  (2004) 299--320.
\newblock \href {http://dx.doi.org/10.1214/154957804100000051}
  {\path{doi:10.1214/154957804100000051}}.

\bibitem{rr04}
G.~Roberts, J.~Rosenthal, General state space markov chains and mcmc
  algorithms, Probab. Surv. 1 (2004) 20--71.
\newblock \href {http://dx.doi.org/10.1214/154957804100000024}
  {\path{doi:10.1214/154957804100000024}}.

\bibitem{gv99}
G.~Latouche, V.~Ramaswami, Introduction to Matrix Analytic Methods in
  Stochastic Modeling, Philadelphia: Society for Industrial Mathematics, 1999.
\newblock \href {http://dx.doi.org/10.1155/S1048953399000362}
  {\path{doi:10.1155/S1048953399000362}}.

\bibitem{s80}
E.~Seneta, Computing the stationary distribution for infinite markov chains,
  Linear Alg. Appl. 34 (1980) 259--267.
\newblock \href {http://dx.doi.org/10.1016/0024-3795(80)90168-8}
  {\path{doi:10.1016/0024-3795(80)90168-8}}.

\bibitem{a03}
S.~Asmussen, Applied Probability and Queues, 2nd Edition, New York:
  Springer-Verlag, 2003.
\newblock \href {http://dx.doi.org/10.1007/b97236} {\path{doi:10.1007/b97236}}.

\bibitem{lzz10}
Y.~Liu, H.~Zhang, Y.~Zhao, Subgeometric ergodicity for continuous-time markov
  chains, J. Math. Anal. Appl. 368 (2010) 178--189.
\newblock \href {http://dx.doi.org/10.1016/j.jmaa.2010.03.019}
  {\path{doi:10.1016/j.jmaa.2010.03.019}}.

\bibitem{Glyn2002}
P.~Glynn, W.~Whitt, Necessary conditions in limit theorems for cumulative
  processes, Stoch. Proc. Appl. 98 (2002) 199--209.
\newblock \href {http://dx.doi.org/10.1016/S0304-4149(01)00146-6}
  {\path{doi:10.1016/S0304-4149(01)00146-6}}.

\bibitem{lz15}
Y.~Liu, Y.~Zhang, Central limit theorems for ergodic continuous-time markov
  chains with applications to single birth processes, Front. Math. China. 10(4)
  (2015) 933--947.
\newblock \href {http://dx.doi.org/10.1007/s11464-015-0488-5}
  {\path{doi:10.1007/s11464-015-0488-5}}.

\bibitem{DS87}
G.~Diana, E.~Seneta, Augmented truncations of infinite stochastic matrices, J.
  Appl. Probab. 24(3) (1987) 600--608.
\newblock \href {http://dx.doi.org/10.1017/S0021900200031338}
  {\path{doi:10.1017/S0021900200031338}}.

\bibitem{cz04}
M.~Chen, Y.~Zhang, Unified representation of formulas for single birth
  processes, Front. Math. China. 9(4) (2014) 761--796.
\newblock \href {http://dx.doi.org/10.1007/s11464-014-0381-7}
  {\path{doi:10.1007/s11464-014-0381-7}}.

\bibitem{jly14}
S.~Jiang, Y.~Liu, S.~Yuan, Poisson's equation for discrete-time single-birth
  processes, Statist. Probab. Lett. 85 (2014) 78--83.
\newblock \href {http://dx.doi.org/10.1016/j.spl.2013.11.008}
  {\path{doi:10.1016/j.spl.2013.11.008}}.

\bibitem{z18}
Y.~Zhang, Criteria on ergodicity and strong ergodicity of single death
  processes, Front. Math. China. 13(5) (2018) 1215--1243.
\newblock \href {http://dx.doi.org/10.1007/s11464-018-0722-z}
  {\path{doi:10.1007/s11464-018-0722-z}}.

\bibitem{lt15}
Y.~Liu, Y.~Tang, Y.~Zhao, Censoring technique and numerical computations of
  invariant distribution for continuous-time markov chains (in chinese), Sci.
  Sin. Math. 45 (2015) 671--682.
\newblock \href {http://dx.doi.org/10.1360/N012015-00074}
  {\path{doi:10.1360/N012015-00074}}.

\bibitem{WZ20}
J.~Wang, Y.~Zhang, Moments of integral-type functionals downward for single
  death processes, Front. Math. China. 15(4) (2020) 749--768.
\newblock \href {http://dx.doi.org/10.1007/s11464-020-0850-0}
  {\path{doi:10.1007/s11464-020-0850-0}}.

\bibitem{c60}
K.~Chung, Markov Chains with Stationary Transition Probabilities, Springer,
  1960.
\newblock \href {http://dx.doi.org/10.1007/978-3-642-49686-8}
  {\path{doi:10.1007/978-3-642-49686-8}}.

\end{thebibliography}
\end{document}